\documentclass[aos-blank]{imsart}
\RequirePackage{amsthm,amsmath,amsfonts,amssymb}
\RequirePackage[numbers]{natbib}
\RequirePackage[colorlinks,citecolor=blue,urlcolor=blue]{hyperref}
\RequirePackage{graphicx}

\usepackage{verbatim}

\startlocaldefs
\theoremstyle{plain}
\newtheorem{axiom}{Axiom}

\newtheorem{theorem}{Theorem}
\newtheorem{lemma}[theorem]{Lemma}
\newtheorem{proposition}[theorem]{Proposition}
\newtheorem{corollary}[theorem]{Corollary}

\newtheorem{example}[theorem]{Example}

\newtheorem{remark}[theorem]{Remark}

\usepackage{enumerate}

\usepackage{tikz}
\usetikzlibrary {shapes}
\usetikzlibrary {arrows}
\usetikzlibrary {positioning}
\usetikzlibrary {calc}
\usetikzlibrary{fit}					
\usetikzlibrary{backgrounds}	

\usepackage{color}

\newcommand{\E}{\mathrm{E}}

\newcommand{\pa}{\mathrm{pa}}

\newcommand{\an}{\mathrm{an}}

\newcommand{\dse}{\,\mbox{$\perp$}\,}

\newcommand{\cd}{\,|\,}

\newcommand{\ci}{\mbox{\protect $\: \perp \hspace{-2.3ex}
\perp$ }}

\newcommand{\notdse}{\nolinebreak{\not\hspace{-1.5mm}\dse}}

\newcommand{\notci}{\nolinebreak{\not\hspace{-1.5mm}\ci}}

\newcommand{\n}[0]{\hspace*{.35em}}

\newcommand{\nn}[0]{\hspace*{.7em}}

\newcommand{\node}{\mbox {\LARGE
{$\mbox{$\circ$}$}}}

\newcommand{\ful}{\mbox{$\, \frac{ \nn \nn \;}{ \nn \nn
}$}}

\newcommand{\fla}{\mbox{$\hspace{.05em} \prec
\!\!\!\!\!\frac{\nn \nn}{\nn}$}}

\newcommand{\fra}{\mbox{$\hspace{.05em} \frac{\nn
\nn}{\nn
}\!\!\!\!\! \succ \! \hspace{.25ex}$}}

\newcommand{\arc}{\mbox{$\hspace{.06em} \prec
\!\!\!\!\!\frac{\nn \nn}{\nn}
\!\!\!\!\!
\succ\! \hspace{.25ex}$}}

\newcommand{\doo}{\mathrm{do}}

\newcommand{\cause}[1]{\mathrm{cause}({#1})}

\newcommand{\dcause}[1]{\mathrm{dcause}({#1})}

\newcommand{\icause}{\mathrm{\iota cause}}

\newcommand{\cc}[1]{\mathrm{cc}({#1})}

\newcommand{\eff}[1]{\mathrm{eff}({#1})}

\newcommand{\neff}[1]{\mathrm{neff}({#1})}

\newcommand{\sco}{\mathrm{sc}}

\newcommand{\ts}[1]{\textcolor{blue}{*** TS: #1 ***}} 

\newcommand{\erk}{\hfill \ensuremath{\Diamond}} 

\newcommand\ns[1]{ \left\{ {#1} \right\} }

 \newcommand\setm[1]{\setminus { \ns{#1} }}

\endlocaldefs

\begin{document}

\begin{frontmatter}
\title{Axiomatization of Interventional Probability Distributions}
\runtitle{Axiomatization of Interventional Distributions}

\begin{aug}
\author[A]
{\fnms{Kayvan}~\snm{Sadeghi}\ead[label=e1]
{k.sadeghi@ucl.ac.uk} \orcid{0000-0001-7314-744X}}
\and
\author[A]
{\fnms{Terry}~\snm{Soo}\ead[label=e2]{math@terrysoo.com} \orcid{0000-0003-1744-1876}}
\address[A]
{Department of Statistical Science,
University College London\printead[presep={,\ }]{e1,e2}}
%
\end{aug}

\begin{abstract}
Causal intervention is an essential tool in causal inference. It is axiomatized under the rules of do-calculus in the case of structure causal models.
We provide simple axiomatizations for families of probability distributions to be different types of interventional distributions.
Our axiomatizations neatly lead to a simple and clear theory of causality that has several advantages: it does not need to make use of any modeling assumptions such as those imposed by structural causal models; it only relies on interventions on single variables; it includes most cases with latent variables and causal cycles; and more importantly, it does not assume the existence of an underlying true causal graph as we do not take it as the primitive object\textemdash in fact, a  causal graph is derived as a by-product of our theory. We show that, under our axiomatizations, the intervened distributions are Markovian to the defined intervened causal graphs, and an observed joint probability distribution is Markovian to the obtained causal graph; these results are consistent with the case of structural causal models, and as a result, the existing theory of causal inference applies.  We also show that a large class of  natural structural causal models satisfy the theory presented here. We note that the aim of this paper is axiomatization of interventional families, which is subtly different from ``causal modeling.''
\end{abstract}

\begin{keyword}[class=MSC]
\kwd[Primary ]{62H22}
\kwd[; secondary ]{62A01}
\end{keyword}

\begin{keyword}
\kwd{ancestral graphs}
\kwd{causal graphs}
\kwd{directed acyclic graphs}
\kwd{do-calculus}
\kwd{interventional distributions}
\kwd{Markov properties}
\kwd{structural causal models}
\end{keyword}

\end{frontmatter}

\section{Introduction}
A popular approach to infer causal relationships is to use the concept of intervention as opposed to observation. For example, as described in \cite{pet17}, it can be observed that there is a correlation between smoking and the colour of the teeth, but no matter how much one whitens somebody's teeth, it would not affect their smoking habits. On the other hand, forcing someone to smoke would affect the colour of their teeth. Hence, smoking has a causal effect on  the colour of the teeth, but not vice versa.

Interventions have generally been embedded in the setting of structural causal models (SCMs), also known as structural equation models \cite{pea88,spioo}. These are a system of assignments for a set of random variables ordered by an associated ``true causal graph,'' which is generally assumed to be unknown. SCMs utilise the theory of graphical (Markov) models, which  are statistical models over graphs with nodes as random variables and edges that indicate some types of conditional dependencies; see \cite{lau96}.

An axiomatic approach to interventions for SCMs, known as Pearl's do-calculus \cite{pea09}, has been developed for identifiability of interventional distributions from the observational ones; see also \cite{hua06,shp06} for some further theoretical developments. There has also been a substantial amount of work on generalizing the concept of intervention from the case of directed acyclic graphs (DAGs) (i.e., Bayesian networks) to more general graphs containing bidirected edges (which indicate the existence of latent variables) (see, e.g., \cite{zha08b}), and directed cycles \cite{bon20}. However, most of these attempts stay within the setting of SCMs or, at least, under the assumption that there exists an  underlying ``true causal graph''  that somehow captures the causal relationships \cite{woo04}.

Interventions in the literature (i.e.\ on SCMs) have been defined to be of different forms; see \cite{kor04,ebe07}. The type of intervention with which we are dealing here is hard in the sense that it destroys all the causes of the intervened variable, and is stochastic in the sense that it replaces the marginal distribution of the intervened variable with a new distribution; although, we will show in Remark \ref{rem:atom} that an atomic (also called surgical) intervention (which forces the variable to have a specific value) can be easily adapted in this setting.

In this paper, without assuming any modeling assumptions such as those given in the setting of SCMs, we give simple conditions for a family of joint distributions $\mathcal{P}_{\doo}= \ns{P_{\doo(1)}, \ldots, P_{\doo(N)}}$ to act as a well-behaved interventional family, so that one can think of $P_{\doo(i)}$ as an interventional distribution on a single variable $X_i$, for each $i\in V=\{1,\ldots,N\}$  in a random vector $X_V$. As apparent from the context, our approach here is aligned with the interventional approach to causality rather than the counterfactual approach.

Here, we are not providing an alternative to the current mainstream setting (sometimes called the Pearlian setting), in the case where it uses the interventional approach, and in most cases, SCMs, and which has led to extensive work on causal learning and estimation. We simply provide theoretical backing for this approach and generalize it beyond SCMs, by providing certain axioms in order to derive as results some of the assumptions that have been used in the literature (such as the existence of the causal graph and the global Markov property w.r.t.\ it).

This paper carries certain important messages: (1) There is no need to take the true causal graph as the primitive object: causal graph(s) can then be formally defined and derived from interventional families (rather than posited). (2) The causal structure (and graph) can be solely derived from the family of interventional distributions; in other words, there is no need for an initial state, i.e., an underlying joint observational distribution $P$ of $X_V$, to be be known for this purpose. However, we provide axioms such that the required consistency between the interventional family and the ``underlying" observational distribution is satisfied when indeed the observational distribution is available, and such that one can measure the causal relationships. (3) To derive the causal structure or graph, (in most situations) one needs to rely only on single interventions once at a time. This is an advantage as much less information is used by only relying on single interventions.  Indeed, there are real world situations in which one would like to consider intervening on several variables simultaneously; we believe a similar theory can be proposed in such cases.

We must emphasize that the work presented here is about axioms that interventional distributions
\emph{should}
satisfy for the purpose of causal reasoning. These axioms should not be confused with a causal model whose goal is to provide ``correct" interpretation of causal relationships and measuring their effects. This difference is quite subtle and could lead to confusion. Similarly, one should distinguish the ``causal graphs" defined and derived here from a graph learned by structure learning from observational and, potentially, interventional data (as in, e.g., \citep{spioo,col14}). The goal here is
\emph{not}
structure learning.

\subsection{Key results} One of our central assumptions, Axiom \ref{prop:trans}, is that cause is transitive; see  \cite{hal00} for a philosophical discussion on the transitivity of the cause.
Under the condition of singleton-transitivity and simple assumptions on  conditional independence structure of $\mathcal{P}_{\doo}$, we show that the causal relations are transitive; see Theorem \ref{prop:intind}.  We provide a definition of causation similar to that in \cite{pet17}. The concept of direct cause is defined  in terms of the conditional independence properties of the interventional family (which is a departure  from the widely-known definition \cite[page 55]{woo04}), and from this we  \emph{define} the intervened causal graph, and using these, we define the causal graph; see Section \ref{sec:caus-gr}; this is a major relaxation of assumptions from the current paradigm  where it is assumed that such a causal graph exists. The obtained causal graph allows bidirected edges and directed cycles without double edges consisting of a bidirected edge and an arrow. We call this family of graphs bowless directed mixed graphs (BDMGs). The generated graph \emph{is} the ``true'' causal graph under the axiomatization, and we show that the definitions related to causal relationships and the graphical notions on the graph are interchangeable (Theorem \ref{thm:exch}).


One of our main theorems is that, under some additional assumptions, namely intersection and composition,  intervened distributions $P_{\doo(i)}$ in the interventional family are Markovian to the defined intervened graphs; see Theorem \ref{thm:markov-int}. We provide additional axioms (Axioms \ref{eq:ax20} and \ref{eq:ax2}) to relate $\mathcal{P}_{\doo}$ to an observed distribution $P$, and call the interventional family (strongly) observable. We show that the underlying distribution $P$ for an observable interventional family is Markovian to the defined causal graph; see Theorem \ref{thm:markov}.  Therefore, the established theory of causality using SCMs, which mainly relies on the Markov property of the joint distribution of the SCM, could be followed from our theory.

We later provide additional axioms (Axioms \ref{eq:ax1} and \ref{eq:ax2n}) for the case of ancestral causal graphs, to define what we call quantifiable interventional families that allow for measuring causal effects. We show that the quantifiable interventional families are strongly observational; see Theorem \ref{thm:biquanttoobint}.

We also compare and contrast our theory with the SCM setting. We show that, for SCMs with certain simple properties (which include transitivity of the cause and are implied by faithfulness), the family of interventions on each node constitutes a strongly observable interventional family, and (even without the transitivity assumption in the case of ancestral graphs) the causal graph generated by the theory presented in this paper is the same as the causal graph associated to the SCM; see Theorem \ref{main-SCM-comp}.

Our theory is based only on intervening on single variables once at a time. We clearly identify cases (which can only be non-maximal and non-ancestral), where this theory may misidentify some direct causes; see Section \ref{sec:intmult}.

\subsection{Related works} To the best of our knowledge, most of the attempts to abstracting intervention or causality based on intervention in general are substantially different from our approach; see for example \cite{ris21} for a category-theory approach, and more recent \cite{par23} for a measure-theoretic approach.

One such attempt is the seminal work of Dawid on the decision theoretic framework for causal inference; for example, see \cite{daw02,daw21}. Our approach share the same concerns and spirit as that of Dawid's in focusing on the interventions rather than counterfactuals, as well as trying to justify the existence of the causal graph rather than assuming it,
as we have heeded Dawid's caution \cite{daw10}.
However, the mechanics of the two works are different.
First of all, we have not provided any statistical model like Dawid does.
In addition, the goal of Dawid's work is mainly to enable ``transfer of probabilistic information from an observational to an interventional setting," whereas, here, our starting point is interventions.
Finally, our approach also covers a much more general class of causal graphs than DAGs, considered by Dawid. We have not included influence diagrams, as proposed by Dawid, in our setting, but believe that it should be possible to derive them together with their conditional-independence constraints using our causal graphs and Markov properties.

A more similar approach to ours is that of
Bareinboim, Brito, and Pearl \cite{bar11}.
Like us the authors start with a family of interventional distributions with interventions defined on single variables. A difference is that they only work on atomic interventions, which requires alternative definitions to conditional independence (which are phrased as invariances)\textemdash we believe this can be adapted to use stochastic intervention and regular conditional independence. Using this, they define the concept of direct cause, which is different, but of similar
to how we define this concept.
A major difference is that, in
their
paper, they provide different notions of compatibility of the interventional family and causal graphs by assuming certain conditions that include the global Markov property\textemdash in our work, we do not assume the global Markov property, and generate a graph directly from the interventional family by using direct cause, and prove the Markov property under certain axioms and conditions.
Another important difference is that
their
work relies on the notion of observed initial state distribution to define interventions\textemdash as mentioned before, we do not need to rely on observational distribution to derive the causal graph.

The mentioned paper is purely on DAGs, but in a more recent work \cite{bar22} the method was generalized to include arcs representing latent variables. The notion of Markov property has been replaced by semi-Markov property to deal with this generalization, but the difference between the two methods remains the same as described for the original paper.
\subsection{Structure of the manuscript}   In the next section, we provide some preliminary material including novel results on the equivalence of the pairwise and global Markov property in our newly defined class of BDMGs, which are used in the subsequent sections.  In Section \ref{sec:intfam}, we define foundational concepts related to interventional families, and provide conditions for the axiom of transitivity of the cause. In Section \ref{sec:causgr}, we define the concepts of direct cause and intervened causal as well as causal graphs, and prove the Markov property of the intervened distribution to the intervened causal graph. In Section \ref{sec:obsint}, we prove the Markov property of the underlying distribution to the causal graph under the observable interventional axioms. In Section \ref{sec:anc}, we specialize the definitions and results for directed ancestral graphs. In Section \ref{sec:qauntint}, we provide additional axioms of quantifiable interventional families for the purpose of measuring causal effects. In Section \ref{sec:scm}, we show how intervention on SCMs fits within the framework of this paper. In Section \ref{sec:intmult}, we provide cases where only single interventions may misidentify certain direct causes. We conclude the paper in Section \ref{sec:summary}. In the \hyperref[appn]{Appendix}, we provide certain proofs of results presented in Section \ref{sec:premmark}, including the proof of equivalence of the pairwise and global Markov property for BDMGs.
\section{Preliminaries}\label{sec:prem}
In this section, we provide the basic concepts of graphical and causal modeling needed in the paper.
\subsection{Conditional measures and  independence}
We will work in the following setting.  Let $V$ be a finite set of size $N$. Let $P$ be a probability measure  on the product measurable space  $\mathcal{X} = \prod_{i \in V} \mathcal{X}_i$.  For $A \subseteq V$, we let $\mathcal{X}_A= \prod_{i \in A} \mathcal{X}_i$ and $P^A$ be the marginal measure of $P$ on $\mathcal{X}_A$ given by
$$P^A(W) = P( W\times \mathcal{X}_{V \setminus A})$$
for all measurable $W \subseteq \mathcal{X}_A$.     We will  use the notation
$i \ci_P j$
to mean that the marginal $P^{\ns{i,j}} = P^{i,j}$ is the product measure $P^i \otimes P^j$ on $\mathcal{X}_i \times \mathcal{X}_j$, so that if $X=(X_1, \ldots, X_N)$ is a random vector
(defined on some probability space $(\Omega, \mathcal{F}, \mathbb{P})$)
taking values on $\mathcal{X}$ with law $P$, then $X_i$ is independent of $X_j$ \cite{daw79,lau96}.

For $x_A \in \mathcal{X}_A$, we let $P(\cdot | x_A) = \mathbb{P}(X \in \cdot | X_A = x_A)$ denote a regular conditional probability  \cite{cha97}  so that in particular, we have the disintegration
$$ P(F) = \int_{\mathcal{X}_A}P(F \cd x_A) d P^A(x_A),$$
for $F \subseteq \mathcal{X}$ measurable.
More generally, consider disjoint subsets $A$, $B$, and $C$ of $V$.  We will often consider the marginal of a conditional measure, and  have a slight abuse of notation that:
$$ P^A( \cdot \cd x_{C}) =  [P( \cdot \cd x_{C})]^A.$$
We write $A\ci_P B\cd C$  to denote that the measure $P^{A, B}(\cdot | x_C)=P^{A \cup B}(\cdot | x_C)$ is a product measure on $\mathcal{X}_A \times \mathcal{X}_B$ for $P^C$-almost all $x_C \in \mathcal{X}_C$, so that if $X$ has law $P$, then $X_A$ is \emph{conditionally independent of $X_B$ given $X_C$}.
Sometimes, we will simply say that $A$ is \emph{conditionally independent} of  $B$ given $C$ in $P$.      In addition, when independence fails,  we write $A\notci_P B\cd C$ and say that $A$ and $B$ are \emph{conditionally dependent} given $C$ in $P$.



\subsection{Structural independence properties of a distribution}
\label{2-prop}

A probability distribution $P$ is always a \emph{semi-graphoid} \cite{pea88}, i.e., it satisfies the four following properties for disjoint subsets $A$, $B$, $C$, and $D$ of $V$:
 \begin{enumerate}[1. ]
    \item $A\ci_P B\cd C$ if and only if $B\ci_P A\cd C$ (\emph{symmetry});
    \item if $A\ci_P B\cup D\cd C$, then $A\ci_P B\cd C$ and $A\ci_P D\cd C$ (\emph{decomposition});
    \item if $A\ci_P B\cup D\cd C$, then $A\ci_P B\cd C\cup D$ and $A\ci_P D\cd C\cup B$ (\emph{weak union});
    \item if $A\ci_P B\cd C\cup D$ and $A\ci_P D\cd C$,
    then $A\ci_P B\cup D\cd C$ (\emph{contraction}).
 \end{enumerate}
Notice that the reverse implication of contraction clearly holds by decomposition and weak union. We so use three different properties of conditional independence that are not always satisfied by probability distributions:
\begin{enumerate}[1. ]
  \setcounter{enumi}{4}
	\item if $A\ci_P B\cd C\cup D$ and $A\ci_P D\cd C\cup B$, then $A\ci_P B\cup D\cd C$ (\emph{intersection});
	\item if $A\ci_P B\cd C$ and $A\ci_P D\cd C$, then $A\ci_P B\cup D\cd C$ (\emph{composition});
	\item
	 if $i\ci_P j\cd C$ and $i\ci_P j\cd C\cup\{k\}$, then $i\ci_P k\cd C$ or $j\ci_P k\cd C$ (\emph{singleton-transitivity}),
\end{enumerate}
where $i$, $j$, and $k$ are single elements. A semi-graphoid distribution that satisfies intersection is called \emph{graphoid}. If the distribution $P$ is a regular multivariate Gaussian distribution, then $P$ is a singleton-transitive compositional graphoid; for example see \cite{stu05} and \cite{pea88}.  If $P$ has strictly positive density, it is always a graphoid; see, for example, Proposition 3.1 in \cite{lau96}.
\begin{remark}
\label{to-get-intersect}
We note that if $P$ has full support over its state space,  then it satisfies the intersection property; for a comprehensive discussion and necessary and sufficient conditions, see \cite{pet15}. \erk
\end{remark}
Finally, we define the concept of ordered stabilities \cite{sad17}. We say that $P$ satisfies \emph{ordered upward-} and \emph{downward-stability} w.r.t.\ an order $\leq$ of $V$ if the following hold:
\begin{enumerate}
    	\item[$\bullet$ \ ] if $i\ci_P j\cd C$,  then $i\ci_P j\cd C\cup\{k\}$ for every $k\in V\setminus\{i,j\}$ such that $i< k$ or $j< k$ (ordered upward-stability);
	\item[$\bullet$ \ ] if $i\ci_P j\cd C$, then $i\ci_P j\cd C\setminus\{k\}$ for every $k\in V\setminus\{i,j\}$ such that $i\nless k$, $j\nless k$, and $\ell\nless k$ for every $\ell\in C\setminus\{k\}$ (ordered downward-stability).
\end{enumerate}

\subsection{Graphs and their properties}
We usually refer to a graph as an ordered pair $G=(V,E)$, where $V$ is the \emph{node} set and $E$ is the \emph{edge} set.
 When nodes $i$ and $j$ are the endpoints of an edge, we call them \emph{adjacent}, and write $i\sim j$, and otherwise $i\nsim j$.

We consider two types of edges:  \emph{arrows} ($i\fra j$) and \emph{bidirected edges} or  \emph{arcs} ($i\arc j$). We do not consider graphs that have simultaneous third type of edge: \emph{undirected edges} or \emph{lines} ($i\ful j$). We only allow for the possibility of multiple edges between nodes when they are arrows in two different directions between $i$ and $j$, i.e., $i\fra j$ and $i\fla j$, which we call \emph{parallel arrows}. This means that we do not allow \emph{bows}, i.e., a multiple edge of arrow and arc, to appear in the graph.

A \emph{subgraph} of a graph $G_1$ is graph $G_2$ such that $V(G_2)\subseteq V(G_1)$ and $E(G_2)\subseteq E(G_1)$ and the assignment of endpoints to edges in $G_2$ is the same as in $G_1$. An \emph{induced subgraph} by nodes $A\subseteq V$ is
the
 subgraph that contains all and only nodes in $A$ and all edges between two nodes in $A$.

A \emph{walk} is a list $\langle v_0,e_1,v_1,\dots,e_k,v_k\rangle$ of nodes and edges such that for $1\leq i\leq k$, the edge $e_i$ has endpoints $v_{i-1}$ and $v_i$. A \emph{path} is a walk with no repeated node or edge. When we define a path, we only write the nodes (and not the edges). A maximal set of nodes in a graph whose members are connected by some paths constitutes a \emph{connected component} of the graph.
A \emph{cycle} is a walk  with no repeated nodes or edges except for $v_0=v_k$.

We call the first and the last nodes \emph{endpoints} of the path and all other nodes \emph{inner nodes}. A path can also be seen as a certain type of connected subgraph of $G$;  a \emph{subpath} of a path $\pi$ is  an induced connected subgraph of $\pi$.
For an arrow $j\fra i$, we say that the arrow is \emph{from} $j$ \emph{to} $i$. We also call $j$ a \emph{parent} of $i$, $i$ a \emph{child} of $j$ and we use the
 notation $\pa(i)$ for the set of all parents of $i$ in the graph.
In the cases of $i\fra j$ or $i\arc j$ we say that there is an arrowhead at $j$ or pointing to $j$. 
A path $\langle i=i_0,i_1,\dots,i_n=j\rangle$ (or a cycle where $i=j$) is \emph{directed} from $i$ to $j$ if all $i_ki_{k+1}$ edges are arrows pointing from $i_k$ to $i_{k+1}$. If there is a directed path from $i$ to $j$, then node $i$ is an \emph{ancestor} of $j$ and $j$ is a \emph{descendant} of $i$. We denote the set of ancestors of $j$ by $\an(j)$;  unlike some authors, we do not allow $j\in\an(j)$. Similarly, we define an ancestor of a set of nodes $A\subset V$ given by $\an(A):= [\bigcup_{j\in A}\an(j)]\setminus A$.  If necessary, we might write $\an_G$ to specify that this is the set of ancestors in $G$.

A \emph{strongly connected component} of a graph is the set of nodes that are mutually ancestors of each other, or it is a single node if that node does not belong to any directed cycle.  It can be observed that nodes of the graph are partitioned into strongly connected components. We denote the members of the strongly connected component containing node $i$ by $\sco(i)$.

A \emph{tripath} is a path with three nodes. The inner node $t$ in each of the  three tripaths
$$i\fra\,
t\fla\,j, \  i\arc\,t\fla\, \ j, \  i\arc\,t\arc\,j$$
is
a \emph{collider} (or a collider node) and the inner node of any other tripath
is a \emph{non-collider} (or a non-collider node) on the tripath or, more generally, on any path of which the tripath is a subpath; i.e.\ a node is a collider if two arrowheads meet. A path is called a \emph{collider} path if all its inner nodes are colliders.

The most general class of graphs that naturally arises from the theory presented in this paper is what we call the \emph{bowless directed mixed graph (BDMG)}, which consists of arrows and arcs, and the only multiple edges are parallel arrows (i.e., they are bowless).  


\emph{Ancestral graphs} \cite{ric02} are graphs with arcs and arrows with no directed cycles and no arcs $ij$ such that $i\in\an(j)$. Acyclic directed mixed graphs (ADMGs) \cite{ric03} are graphs with arcs and arrows with no directed cycles. In other words, BDMGs unify ADMGs without bows and directed cycles. BDMGs also trivially contain the class of \emph{directed ancestral graphs}, i.e., \emph{ancestral graphs} \cite{ric02} without lines. 
These all also contain \emph{directed acyclic graphs} (DAGs) \cite{kii84}, which are graphs with only arrows and no directed cycles.

The class of BDMGs is a subclass of \emph{directed mixed graphs}, introduced in \cite{bon20}, which is a very general class of graphs with arrows and arcs that allow for directed cycles. Later on, we will use some definitions and results originally defined for directed mixed graphs.
\subsection{Markov properties}\label{sec:premmark}
\label{sigma-sep}
In this paper, we will need global and pairwise Markov properties for BDMGs.
%
In order to introduce the global Markov property, we need to define the concept of $\sigma$-separation for directed mixed graphs (in fact, originally defined for the larger class of \emph{directed graphs with hyperedges}
in
\cite{forre2017markov}).

A path $\pi=\langle i=i_0,i_1,\cdots,i_n=j\rangle$ is said to be \emph{$\sigma$-connecting given $C$}, which is disjoint from $i,j$, if all its collider nodes are in $C\cup\an(C)$; and all its non-collider nodes $i_r$ are either outside $C$, or if there is an arrowhead at $i_{r-1}$, then $i_{r-1}\in\sco(i_r)$ and if there is an arrowhead at $i_{r+1}$ on $\pi$ then $i_{r+1}\in\sco(i_r)$.
For disjoint subsets $A,B,C$ of $V$, we say that $A$ and $B$ are \emph{$\sigma$-separated given $C$}, and write $A\dse_{\sigma}B\cd C$, if there are no $\sigma$-connecting paths between $A$ and $B$ given $C$.

In the case where there are no directed cycles in the graph, $\sigma$-separation reduces to the $m$-separation of \cite{ric02}; recall that \emph{$\pi$ is $m$-connecting given $C$} if all its collider nodes are in $C\cup\an(C)$; and all its non-collider nodes are outside $C$. In addition, if there are no arcs in the graph, i.e., the graph is a DAG, it reduces to the well-known $d$-separation \cite{pea88}.

We call two graphs \emph{ Markov equivalent} if they induce the same set of conditional separations.

A probability distribution $P$ defined over $V$ satisfies the \emph{global Markov property} w.r.t.\ a bowless directed mixed graph $G$, or is simply \emph{Markovian} to $G$, if for  disjoint subsets $A$, $B$, and $C$ of $V$, we have
$$A\dse_{\sigma} B\cd C \implies A\ci B\cd C.$$
If $G$ is an ancestral graph (or a DAG), then $\dse_{\sigma}$ will be replaced by $\dse_m$ (or $\dse_d$) in the definition of the global Markov property.

If, in addition to the global Markov property, the other direction of the implication holds, i.e.,  $A\dse_{\sigma} B\cd C \iff A\ci B\cd C$,  then we say that  $P$ and $G$ are \emph{faithful}. A weaker condition of \emph{adjacency-faithfulness} \cite{ram06,zha08} states that for every edge between $k$ and $j$ in $G$, there are no independence statements $k\ci_P j \cd C$ for any $C$.

We now define that a distribution $P$ satisfies the \emph{pairwise Markov property} \eqref{P} w.r.t.\ a bowless directed mixed graph $G$, if for every pair of non-adjacent nodes $i,j$ in $G$, we have
\begin{equation}
\tag{\text{PMP}}
\label{P}
i \ci_P j\cd \an(\{i,j\}).
\end{equation}
This is the same wording as that of the pairwise Markov property for the subclass of ancestral graphs; see \cite{sadl16}.

We prove the equivalence of the pairwise and global Markov properties, which shall be used later for causal graphs and Theorem \ref{thm:markov}.   The proofs are presented in the \hyperref[appn]{Appendix} as they are not the main focus of this manuscript.


\begin{theorem}\label{prop:pairMark}
Let $G$ be a BDMG, and $P$ satisfy the intersection and composition properties.  If $P$ satisfies the pairwise Markov property \eqref{P} with respect to $G$, then $P$ is Markovian to $G$.
\end{theorem}

We also define the converse of the pairwise Markov property. We say  that $P$ satisfies the \emph{converse pairwise Markov property} w.r.t.\ $G$  if an edge between $i$ and $j$ in $G$ implies that
$$i \notci_P j\cd \an(\{i,j\}).$$
Notice that faithfulness and adjacency-faithfulness of $P_{\mathcal{C}}$ and $G_{\mathcal{C}}$ imply the converse pairwise Markov property; see \cite{sad22}.

A graph is called \emph{maximal} if the absence  of an edge between $i$ and $j$ corresponds to a conditional separation statement for $i$ and $j$, i.e.\ there exists for some $C$ a statement of form $i\dse j\cd C$. Notice from the definition of $\sigma$-separation that graphs with \emph{chordless} directed cycles (i.e., having two non-adjacent nodes in a cycle) are not maximal.

We have the following corresponding converse to Theorem \ref{prop:pairMark}.

\begin{proposition}
\label{prop:Markpair}
Let $G$ be a maximal BDMG. If $P$ is Markovian to $G$,  then $P$ satisfies the pairwise Markov property \eqref{P} with respect to $G$.
\end{proposition}

We call a non-adjacent pair of nodes which cannot be $\sigma$-separated in a non-maximal graph, regardless of what to condition on,
an
\emph{inseparable} pair, and  a non-adjacent pair of nodes which can be $\sigma$-separated in a maximal or non-maximal graph a \emph{separable} pair.

For BDMGs, we define a \emph{primitive inducing path (PIP)}  to be a  path $\langle i,q_1,\cdots,q_r,j\rangle$ (with at least $3$ nodes) between $i$ and $j$, where
\begin{enumerate}[(i) ]
    \item all edges $q_mq_{m+1}$ are either arcs or an arrow where $q_m\in\sco(q_{m+1})$ except for the first and last edges, which may be $i\fra q_1$ or $q_r\fla j$ (without being in the same connected component);
    \item for all inner nodes, we have  $q_m\in\an(\{i,j\})$, i.e., they are in ancestors of $i$ or $j$.
\end{enumerate}
PIPs were originally defined for the case of ancestral graphs (where they were allowed to be an edge) \cite{ric02}.
We show in the \hyperref[appn]{Appendix} the result below:
\begin{proposition}\label{prop:nonmaxp}
In a BDMG,  inseparable pairs are connected by PIPs.
\end{proposition}

In addition, notice from Proposition~\ref{prop:Markpair} that if $P$ is Markovian to $G$, then a pair $i,j$ being a separable pair is equivalent to the separation $i\dse_{\sigma}j\cd \an{\{i,j\}}$.
%

We say that a  graph $G=(V,E)$ admits a \emph{valid order} $\leq$ if for nodes $i$ and $j$ of $G$, $i\fra j$ implies that $i >j$; and $i\arc j$ implies that $i$ and $j$ are incomparable. Notice that this specifies the partial order via its cover relations. In fact, this order can is used as the order w.r.t.\ which ordered upward- and downward-stability hold for graph separations; see \cite{sad17}.

Finally, for ancestral graphs, the $m$-separation holds for every set ``between'' parents and ancestors. We will use this fact later:
\begin{lemma}
\label{lem:squeeze}
For an ancestral graph and  for separable nodes $i,j$, we have
$$i\dse_m j\cd A,$$
for every $A$ such that $\pa(\{i,j\})\subseteq A\subseteq \an(\{i,j\})$.
\end{lemma}
\begin{proof}
The set $\pa(i)$ separating $i$ and $j$ follows from the fact that it is the \emph{Markov blanket}; see \cite{ric03}. The fact that $m$-separation satisfies ordered upward-stability w.r.t.\ orderings associated to an ancestral graphs \cite{sad17} implies that $i\dse_m j\cd \pa(\{i,j\})$, and consequently, $i\dse_m j\cd A$.
\end{proof}
\subsection{Structural causal models} The theory in this paper does not use or assume \emph{structural causal models} (SCMs) (also known as the \emph{structural equation models}) \cite{pea09,spioo}. We define SCMs here as they are an interesting special case of our theory, for which intervention could be easily conceptualized.

Here, we define SCMs for the class of BDMGs as a simplified version of SCMs defined for directed (mixed) graphs in \cite{bon20}. Consider a graph $G$ with $N$ nodes, which in the context of causal inference is often referred to as the ``true causal graph.'' A structural causal model $\mathfrak{C}$  associated with $G$ is defined as a collection of $N$ equations
\begin{equation*}\label{eq:scm}
X_i=\phi_i(X_{\pa_G(i)},\epsilon_i), \nn i\in\{1,\dots, N\},
\end{equation*}
where $\pa_G(i)$ is defined on $G$ and $\epsilon_i$ might be called \emph{noises}; for any subsets $A$ and $B$, we require that $\epsilon_A \ci \epsilon_B$  if and only if, in $G$, there is no arc between any node in $A$ and any node in $B$. In this paper, we usually refer to an SCM as $\mathcal{C}$ and its joint distribution as $P_{\mathcal{C}}$.

In the more widely-used case where $G$ is a DAG, all the $\epsilon_i$ are jointly independent.  For both mathematical and causal discussions on SCMs with DAGs, see \cite{pet17}. When directed cycles are existent, some solvability conditions are required in order for the theory of SCMs to work properly; for this and for more general discussion, see \cite{bon20}.

\emph{Standard interventions} are defined quite naturally when functional equations are specified, as in the case of  SCMs: By intervening on $X_i$ we replace the equation associated to $X_i$ by $X_i=\tilde{X}_i$, where $\tilde{X}_i$ is independent of all other noises; it is not necessary that $\tilde{X}_i$ has the same distribution as $X_i$.   We are concerned with a similar type of intervention in this paper -- this  is a special case of the so-called \emph{stochastic intervention} \cite{kor04}, where some parental set of $X_i$ might still exist after intervention on $X_i$; see also \cite{pet17}. A more special type of intervention is called \emph{perfect intervention} (or \emph{surgical intervention}), where it puts a point mass on a real value $a$ -- this is the original idea of \emph{do-calculus} \cite{pea09}, and is often denoted by $\doo(X_i=a)$.

An important result for SCMs, which facilitates causal inference using them, is that the joint distribution of an SCM is Markovian to its associated graph; see \cite{ver88,pea09} for the case of DAGs, \cite{sad22} for \emph{directed ancestral graphs}, and \cite{bon20} for directed mixed graphs.
\section{Interventional family of distributions}\label{sec:intfam}
\subsection{Interventional families and the cause}
Again, let  $V$ be a finite set of size $N$. 
 Consider a family of distributions $\mathcal{P}_{\doo}=\{P_{\doo(i)}\}_{i\in V}$, where each  $P_{\doo(i)}$ is defined over the same state space  $\mathcal{X} = \prod_{i \in V} \mathcal{X}_i$. We refer to $\mathcal{P}_{\doo}$ as an \emph{interventional family} (of distributions). For $(\tilde{X}_j)_{j \in V}$  a random vector with  distribution $P_{\doo(i)}$, we think of $P_{\doo(i)}$ as the interventional distribution after intervening on some variable $X_i$.

For $k\in V$, we define the set of the \emph{causes} of $k$ as
$${\rm{cause}}^{\mathcal{P}_{\doo}}(k)=\cause{k}:=\{i:\n i\neq k,\n i\notci_{P_{\doo(i)}}k\};$$
we  rarely, simultaneously, have to consider two different interventional families at once.
   Thus, if $i \in \cause{k}$, then, for $(\tilde{X}_j)_{j \in V}$, we have that $\tilde{X}_k$ is dependent on $\tilde{X}_i$. Notice that, by convention, $k\notin \cause{k}$. 
   For a subset $A\subseteq V$, we define $\cause{A}=(\bigcup_{k\in A} \cause{k})\setminus A$.

%

%
%
%

The definition of the cause  is identical to what is know in the literature as the ``existence of the total causal effect" (see, e.g., Definition 6.12 in \cite{pet17}). Its combination with $\mathcal{P}_{\doo}$ meets the intuition behind cause and intervention: after intervention on a variable $X_i$, it is dependent on a variable $X_j$ if and only if it is a cause of that variable. 


The above setting is compatible with the well-known intervention for SCMs; 
 for a comprehensive discussion on this, see Section \ref{sec:scm}. 
 We will often illustrate our theory with simple examples of SCMs and standard intervention on a single node.

We can also define the set of \emph{effects} of $i$
denoted by $\eff{i}$
by $i\in\cause{k}\iff k\in\eff{i}$.   We take note of the following useful fact.

\begin{remark}\label{rem:eff}
From the definition of cause, we have
$\eff{i}=\{k: i\neq k,\n i \notci_{P_{\doo(i)}}k\}.$
\erk
\end{remark}

\begin{remark}[Interventional families with the same cause]
\label{similar-interv-remark}
Interventions families are defined for one joint distribution per intervention. This can be considered an advantage as not all interventions have to follow a single causal graph. Here, we provide an immediate condition for interventional families to define the same set of causes. After causal graphs have defined, in Section  \ref{sec:cong}, we provide conditions for families of distributions that lead to the same causal graph.

Consider the interventional families $\mathcal{P}_{\doo}=\{P_{\doo(j)}\}_{j\in V}$ and $\mathcal{Q}_{\doo}=\{Q_{\doo(j)}\}_{j\in V}$ over the same state space $\mathcal{X}$.
Causes and effects depend on the interventional family, and, by Remark \ref{rem:eff}, it follows that, for all $i,k \in V$,   
$$\Big[i\notci_{P_{\doo(i)}}k\iff i \notci_{Q_{\doo(i)}}k\Big]  \text{ if and only if }
{\rm{cause}}^{\mathcal{P}_{\doo}}(\ell) = {\rm{cause}}^{\mathcal{Q}_{\doo}}(\ell),$$
for all $\ell\in V$, in which case we say they have the \emph{same causes}.
In the case of SCMs, under standard interventions, the causes are usually invariant with respect to the choice of distribution, except for technical counterexamples; see Remark \ref{peters-thing} and Example \ref{supports-counter}.

 Consider a fixed $i \in V$, and suppose that measures $P_i$ and $Q_i$ have the same null sets.  Notice that the equality
 \begin{equation*}
\label{similar-int}
 Q_{\doo(i)}(\cdot \cd x_i) =  P_{\doo(i)}(\cdot \cd x_i)
 \end{equation*}
for almost every $x_i \in \mathcal{X}_i$, is sufficient for  the effects of $i$ to be same in both families.  Moreover, with the disintegration
$$ dP_{\doo(i)}(x) =  dP_{\doo(i)}(x_{V \setminus \ns{i}} \cd x_i) dP_{\doo(i)}^i(x_i)$$
we see that effects of $i$ depend only on the corresponding conditional distribution, and are invariant under the marginal distributions on $i$ with the same null sets.
\erk
\end{remark}

\begin{remark}[Atomic interventions]\label{rem:atom}
Observe that the dependence
\begin{equation*}
i\notci_{P_{\doo(i)}} k 
\end{equation*}
is equivalent to the existence of  disjoint  measurable subsets $W^{*}, W^{**} \subset \mathcal{X}_i$ of positive measures under  $P_{\doo(i)}^{i}$ satisfying the inequality
\begin{equation}
\label{atomic}
P_{\doo(i)}^{k} (\cdot \cd x_i^*)\neq P_{\doo(i)}^{k}(\cdot \cd x_i^{**}),
\end{equation}
for all $ x_i^* \in W^{*}$ and all $x_i^{**} \in W^{**}$.
Since  $P_{\doo(i)} (\cdot \cd x_i^*)$ and $P_{\doo(i)}(\cdot \cd x_i^{**})$ are probability measures on $\mathcal{X}_{V \setminus \ns{i}}$, they can be thought of as atomic interventions on $i$, where the values at $i$ are fixed, at $x_i^*$ and $x_i^{**}$, respectively.  Thus inequality \eqref{atomic} has the interpretation that $i$ is a cause of $k$ if and only if there exists atomic interventions that witness an effect on $k$; that is, as a function of $x_i$, the conditional probability, $P_{\doo(i)}^{k} (\cdot \cd x_i)$, is non-constant.

From Remark \ref{similar-interv-remark}, without loss of generality,  given  atomic interventions, $A_{\doo(x_i)}$, which are measures on $\mathcal{X}_{V \setminus \ns{i}}$ indexed by $x_i \in \mathcal{X}_i$, we can extend these to an intervention, defined on the complete space, $\mathcal{X}$,  via  the disintegration
$$d{P}_{\doo(i)}(x)    :=  dA_{\doo(x_i)} (x_{V \setminus \ns{i}})dR(x_i),$$
where $R$ is a suitably chosen probability measure on $\mathcal{X}_i$  Specifically, in the case where $\mathcal{X}_i$ is finite, $R$ can be taken to be a uniform measure on $\mathcal{X}_i$.   \erk
\end{remark}

We call a subset $S\subseteq V$ a \emph{causal cycle} if for every $i,k\in S$, we have $k\in\cause{i}\cap \eff{i}$. We write $\cc{i}$ to denote the causal cycle containing $i$.

Under the composition property,  we have the following independence.
 Let $\neff{i}:=V\setminus\eff{i}$ denote the subset of $V$ that contains members that are \emph{not} an effect of $i$.
\begin{proposition}[Non-effects under composition]
\label{prop:intindcomp}
Let $\mathcal{P}_{\doo}$ be  a family of distributions.  If $P_{\doo(i)}$ satisfies the composition property, then
$$i\ci_{P_{\doo(i)}} \neff{i}.$$
\end{proposition}
\begin{proof}
Let $k\in\neff{i}$. Since $i$ is not a cause of $k$, by definition, $i\ci_{P_{\doo(i)}}k$; iteratively applying  the composition property, we obtain the desired result.
\end{proof}

\subsection{Transitive interventional families}

We now say that $\mathcal{P}_{\doo}$ is a \emph{transitive interventional} family if the following axiom holds.
\begin{axiom}[Transitivity of cause]
\label{prop:trans}
For distinct $i,j,k\in V$, if $i\in\cause{j}$ and $j\in\cause{k}$, then $i\in\cause{k}$.
\end{axiom}

Notice also that $P_{\doo(i)}$ in transitive interventional families are restricted: Axiom \ref{prop:trans}  places constraints between different $P_{\doo(i)}$ since $\cause{k}$ depends on all $P_{\doo(i)}$.

Under singleton-transitivity, we have sufficient conditions for  Axiom \ref{prop:trans} to hold. Notice that that these conditions are not satisfied in general, even in the case of an SCM with standard interventions.
\begin{theorem}[Transitivity of cause under singleton-transitivity]
\label{prop:intind}
Let $\mathcal{P}_{\doo}$ be an interventional family whose members $P_{\doo(i)}$ satisfy singleton-transitivity. Assume, for distinct $i,j,k\in V$ such that $i\notin\cause{k}$ and $j \in \cause{k}$, we have:
\begin{enumerate}[(a) ]
  \item
  \label{indep-all}
  $i\ci_{P_{\doo(i)}} k\cd j$; and
  \item
  \label{indep-trans}
  $j\notci_{P_{\doo(i)}} k.$
\end{enumerate}
Then $\mathcal{P}_{\doo}$ is a transitive interventional family.
\end{theorem}
\begin{proof}
Towards a contradiction, assume that $i \in\cause{j}$ and $j\in\cause{k}$, but $i\notin\cause{k}$.  Since  $j\in\cause{k}$,  we have, by \eqref{indep-trans}, that $j\notci_{P_{\doo(i)}} k$.  Also, since $i \in \cause{j}$, we have $i\notci_{P_{\doo(i)}}j$.

The dependencies $j\notci_{P_{\doo(i)}} k$ and $i\notci_{P_{\doo(i)}}j$, along with  singleton-transitivity, in its contrapositive form, imply that $i\notci_{P_{\doo(i)}} k$ or $i\notci_{P_{\doo(i)}} k\cd j$; however, the former is ruled out by assumption, and the latter contradicts \eqref{indep-all}.
%
\end{proof}

\begin{proposition}
\label{coro:order}
Transitivity of the cause induces  a strict preordering $\lesssim$ on $V$ by
\begin{equation}
\left\{
                                  \begin{array}{l}
                                    i<k \iff k\in\cause{i} \text{ and } k\notin\cause{i}; \\
                                    i\sim k \iff k\in\cause{i}\cap \eff{i}.
                                  \end{array}
\right.
\end{equation}
\end{proposition}
\begin{proof}
Irreflexivity is implied by the convention that $i\notin \cause{i}$. Transitivity is Axiom \ref{prop:trans}.
\end{proof}

\begin{corollary}
\label{composition-wu}
Let $\mathcal{P}_{\doo}$ be  a transitive interventional family that satisfies the composition property. If  $i\notin\cause{k}$, then
$$i\ci_{P_{\doo(i)}}k\cd\cause{k}.$$
\end{corollary}
\begin{proof}
Transitivity implies that $i$ is not a cause of members of $\cause{k}$, so that  $\cause{k}\subseteq \neff{i}$.  Thus Proposition \ref{prop:intindcomp}, and the weak union property yield the result.
\end{proof}
In the next example, we show that we cannot drop the singleton transitivity assumption in Theorem \ref{prop:intind}.

\begin{example}[Failure of transitivity without singleton transitivity]  Suppose $1$ is the  cause of $2$, and $2$ is a   cause of $3$.  It may not be the case that $1$ is a cause of $3$.  Consider the SCM, with $X_1 \fra X_2 \fra X_3$, where $X_1$ is Bernoulli $p \in (0,1)$,  conditional on $X_1=x_1$, we sample a Poisson random variable $N=n$ with mean $x_1+1$, and then we sample  $X_2=(X_2^0,\ldots, X_{2}^n)$ as an i.i.d.\ sequence of $n+1$ Bernoulli random variable(s) with parameter $\tfrac{1}{2}$, and finally set $X_3 := X_2^0$.  The first random variable $X_1$ is independent of the final result $X_3$.  The standard interventions where we simply substitute a distributional copy of $X_i$ for each $i$ gives that $1$ is clearly the  cause of $2$, and $2$ a cause of $3$.  However, since $P_{\doo(1)} =P$, singleton-transitivity fails:  $1 \ci_P 3$ and $1 \ci_P 3 \cd 2$, but we have neither $1 \ci_P 2$ nor $3 \ci_P 2$.
\erk
\end{example}

\section{Causal graphs and intervened Markov properties}\label{sec:multsec}
\label{sec:causgr}
\subsection{Intervened and direct cause}
Notice that only knowing the causal ordering cannot yield a graph, since, for example, for $i\fra j\fra k$, there is no way to distinguish the two graphs corresponding to whether an additional $i\fra k$ exists in the graph or not. For this reason, we need to define the concept of direct cause.

In order to define direct and intervened cause in the general case, we need an iterative procedure:
\begin{enumerate}[1. ]
  \item For each $i, k \in V$,  start  with $\dcause{k}:=\cause{k}$, $\icause_i(k):=\cause{k}$, and $S(\mathcal{P}_{\doo})$ an empty graph with node set $V$;
  \item Redefine $$\dcause{k}:=\{i:\n i\in\dcause{k},\n i\notci_{P_{\doo(i)}} k\cd \icause_i(k)\setminus\{i\}\};$$
  \item Generate  $S(\mathcal{P}_{\doo})$ by setting arrows from $i$ to $k$, i.e., $i\fra k$ if $i\in \dcause{k}$;
  \item Generate the graph $S_i(\mathcal{P}_{\doo})$ by removing all arrows pointing to $i$ from $S(\mathcal{P}_{\doo})$;
  \item Redefine $\icause_i(k):=\an_{S_i(\mathcal{P}_{\doo})}(k)$;
  \item If  $S(\mathcal{P}_{\doo})$ is modified by Step 3, then go to Step 2; otherwise, output $\dcause{k}$, $\icause_i(k)$, $S(\mathcal{P}_{\doo})$, and $S_i(\mathcal{P}_{\doo})$.
\end{enumerate}

We call $\dcause{k}$ the set of the \emph{direct causes} of $k$, and $\icause_i(k)$ the set of \emph{intervened causes} of $k$ after intervention on $i$. We also call $S(\mathcal{P}_{\doo})$ the \emph{causal structure}, and $S_i(\mathcal{P}_{\doo})$ the \emph{$i$-intervened causal structure}.

Notice that, since in the iteration $\dcause{k}$ is getting smaller, the procedure will stop.

Notice also that, by convention, $k\notin \dcause{k}\cup\icause_i(k)$, for any $i \in V$.
We also let $\dcause{A}=(\bigcup_{k\in A} \dcause{k})\setminus A$, and $\icause_i(A)=(\bigcup_{k\in A} \icause_i(k))\setminus A$.

As seen by definition, ``cause'' is a \emph{universal} concept: no matter how large the system of random variables is, as long as it contains the two investigated random variables, the marginal dependence of those variables stays intact. On the other hand, ``direct cause'' depends on the system of variables in which the two investigated variables lie.

 The below example shows why we need $\icause_i(k)$ as opposed to $\cause{k}$ in the definition of $\dcause{k}$; see also Section \ref{sec:multint}.
\begin{example}\label{ex:dcauseit}
Let the graph of Figure \ref{fig:dcause}, below,  be the graph associated to an SCM with standard interventions. Under faithfulness, notice that,  $i\notci_{P_{\doo(i)}} k\cd \cause{k}\setminus\{i\}$, where $\cause{k}\setminus\{i\}=\{j,\ell\}$. However, $i$ is clearly not a direct cause of $k$. On the other hand, $\icause_i(k)\setminus\{i\}=\{j\}$, and $i\ci_{P_{\doo(i)}} k\cd \icause_i(k)\setminus\{i\}$.

In the iterative procedure above, in the first round, there will be an arrow from $i$ to $k$ in  $S(\mathcal{P}_{\doo})$. In the second round, $i$ will be removed.
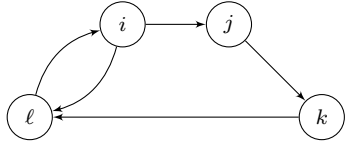
\begin{figure}[htb]
\centering
\begin{tikzpicture}[node distance = 8mm and 8mm, minimum width = 6mm]
    \begin{scope}
      \tikzstyle{every node} = [shape = circle,
      font = \scriptsize,
      minimum height = 6mm,
      inner sep = 2pt,
      draw = black,
      fill = white,
      anchor = center],
      text centered]
      \node(i) at (0,0) {$i$};
      \node(j) [right = of i] {$j$};
      \node(k) [below right = of j] {$k$};
      \node(l) [below left = of i] {$\ell$};
    \end{scope}
		
    \begin{scope}[->, > = latex']
    \draw (i) -- (j);
      \draw (j) -- (k);
    \draw (k) -- (l);
    \draw (l) to [bend left]  (i);
    \draw (i) to [bend left]  (l);
    \end{scope}

    \end{tikzpicture}
		\caption{{\footnotesize A graph for which detecting the direct cause requires an iterative procedure.}}\label{fig:dcause}
		\end{figure} \erk
\end{example}

\begin{remark}\label{rem:dcause}
For causal graphs (as defined in the next subsection) that are ancestral, $\dcause{k}$ can simply be defined by
$$i\notci_{P_{\doo(i)}} k\cd \cause{k}\setminus\{i\}$$
rather than using the iterative procedure and $\icause_i(k)$ in the conditioning set; see Section \ref{sec:anc} for the equivalence of the two methods under certain conditions. \erk
\end{remark}

\begin{remark}
As it is seen in Section \ref{sec:intmult}, there may still be arrows generated here that arguably should not be considered direct causes in the case of non-maximal non-ancestral graphs. We study and identify these cases and offer adjustments in that section. \erk
\end{remark}
The next proposition can be thought of as a proxy for the pairwise Markov property, and will be used in our proofs of the Markov property for our casual  graphs.

\begin{proposition}\label{prop:axcon}
Let $\mathcal{P}_{\doo}$ be a transitive interventional family, and $P_{\doo(i)}$ satisfy the composition property.
For distinct $i,k\in V$,  if $i\notin\dcause{k}$, then
  $$i\ci_{P_{\doo(i)}} k\cd \icause_i(k)\setminus\{i\}.$$
\end{proposition}
\begin{proof}
If $i\in\cause{k}$, then the result follows directly from the definition of $\dcause{k}$.
If $i\notin\cause{k}$, then by Corollary \ref{composition-wu}, we have $i\ci_{P_{\doo(i)}}k\cd\cause{k} \setminus{\ns{i}}$; moreover,
observe that transitivity also implies that  $\icause_i(k)=\cause{k}$.
\end{proof}
\subsection{Intervened causal and causal graphs}\label{sec:caus-gr}

We are now ready to define a graph that demonstrates the causal relationships by capturing the direct causes as well as non-causal dependencies due to latent variables.

Given an interventional family $\mathcal{P}_{\doo}$, and $i\in V$, we define the \emph{$i$-intervened graph}, denoted by $G_i(\mathcal{P}_{\doo})$ to be the $i$-intervened causal structure, where, in addition,  for each pair of non-adjacent nodes $j,k \in V$ that are distinct from $i$, we place an arc between $j$ and $k$, i.e., $j\arc k$, if
\begin{equation}
\label{put-arc}
j\notci_{P_{\doo(i)}} k\cd \icause_i(\{j,k\}),
\end{equation}

Thus with \eqref{put-arc}, we put an arc if \emph{one} of the interventions $P_{\doo(i)}$ suggests the presence of a latent variable.

We also define the \emph{causal graph}, denoted by $G(\mathcal{P}_{\doo})$, to be the  causal structure, where, in addition,  for each pair of nodes $j,k$ that are not adjacent by an arrow, we place an arc between them if the $jk$-arc exists in $G_i(\mathcal{P}_{\doo})$ for \emph{every} $i\in V$ that is distinct from $j$ and $k$.

\begin{remark}
We note that $G_i(\mathcal{P}_{\doo})$ does not contain arrows pointing to $i$ in $G(\mathcal{P}_{\doo})$ and all arcs with $i$ as an endpoint in $G(\mathcal{P}_{\doo})$.

We also note that the existence of the $jk$-arc in two different intervened graphs $G_i(\mathcal{P}_{\doo})$ and $G_\ell(\mathcal{P}_{\doo})$ may not coincide when there is a PIP as the below example shows. \erk
\end{remark}
\begin{example}
Consider an SCM with the graph presented in Figure \ref{fig:nonmax}, below,  with standard intervention. Assume that the joint distribution of the SCM is faithful to this graph. It is easy to observe that $j\notci_{P_{\doo(i)}} k\cd \icause_i(\{j,k\})$, but $j\ci_{P_{\doo(\ell)}} k\cd \icause_i(\{j,k\})$. This implies that  there is a $jk$-arc in $G_i(\mathcal{P}_{\doo})$, but not in $G_\ell(\mathcal{P}_{\doo})$.
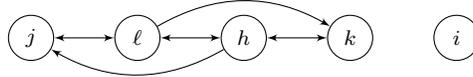
\begin{figure}[htb]
\centering
\begin{tikzpicture}[node distance = 8mm and 8mm, minimum width = 6mm]
    \begin{scope}
      \tikzstyle{every node} = [shape = circle,
      font = \scriptsize,
      minimum height = 6mm,
      inner sep = 2pt,
      draw = black,
      fill = white,
      anchor = center],
      text centered]
      \node(l) at (0,0) {$\ell$};
      \node(h) [right = of l] {$h$};
      \node(k) [right = of h] {$k$};
      \node(j) [left = of l] {$j$};
      \node(i) [right = of k] {$i$};
    \end{scope}
		
    \begin{scope}[->, > = latex']
    \draw (l) to [bend left]  (k);
    \draw (h) to [bend left]  (j);
    \end{scope}
    \begin{scope}[<->, > = latex']
    \draw (j) -- (l);
    \draw (l) -- (h);
    \draw (h) -- (k);
    \end{scope}

    \end{tikzpicture}
		\caption{{\footnotesize A (non-maximal) graph associated to an SCM}}\label{fig:nonmax}
		\end{figure} \erk
\end{example}
For an interventional family $\mathcal{P}_{\doo}$, our definitions now allow us to use the notions $\pa_{G(\mathcal{P}_{\doo})}(k)$ and  $\dcause{k}$  interchangeably. For a transitive $\mathcal{P}_{\doo}$,  we have that, moreover, other causal terminologies on $\mathcal{P}_{\doo}$ and the graph terminologies on $G(\mathcal{P}_{\doo})$ can be used interchangeably:
\begin{theorem}[Interchangeable  terminology]\label{thm:exch}
For a transitive interventional $\mathcal{P}_{\doo}$ where $P_{\doo(i)}$ satisfy the composition property, we have the following:

\begin{enumerate}[(i) ]
  \item
  $i\in\an_{G(\mathcal{P}_{\doo})}(k) \Longleftrightarrow i\in \cause{k} \text{ in } \mathcal{P}_{\doo}.$
  \item
  $i\in\sco_{G(\mathcal{P}_{\doo})}(k) \Longleftrightarrow k\in\cc{i} \text{ in } \mathcal{P}_{\doo}$. \hfill
\end{enumerate}
\end{theorem}
\newpage
\begin{proof}
\hspace{1 cm}
\begin{enumerate}[(i) ]
\item
 The direction $i\in\an_{G(\mathcal{P}_{\doo})}(k) \Rightarrow i\in \cause{k}$ is implied by the  transitivity of ancestors and Axiom \ref{prop:trans}.

Assume that $i\in \cause{k}$.  We prove the result by induction on the cardinality of  $\icause_i(k)\cap\eff{i}$.  For the basis, if   $|\icause_i(k)\cap\eff{i}|=0$, then by the definition of $\eff{i}$, for each $\ell \in \icause_i(k) \setm{i}$, we must have $i\ci_{P_{\doo(i)}} \ell$; furthermore, composition implies $$i\ci_{P_{\doo(i)}}\icause_i(k) \setm{i}.$$ Towards a contradiction, suppose that  $i\notin\dcause{k}$, so that
$$i\ci_{P_{\doo(i)}}k \cd\icause_i(k) \setm{i}.$$
Contraction and decomposition imply $i\ci_{P_{\doo(i)}} k$, so that we have $i \not \in \cause{k}$,  contrary to our original assumption.

For the inductive step,  assume that $|\icause_i(k)\cap\eff{i}|=r \geq 1$. Choose an $\ell\in\icause_i(k)\cap\eff{i}$. The transitivity of the cause, assumed in Axiom \ref{prop:trans}, implies that $\icause_i(\ell)\subset\icause_i(k)$. This implies that $|\icause_i(\ell)\cap\eff{i}|<r$, which, by inductive hypothesis, implies there is a directed path from $i$ to $\ell$.

Similarly,  $\eff{\ell}\cap\icause_i(k)\subset\eff{i}\cap\icause_i(k)$, so that $|\eff{\ell}\cap\icause_i(k)|<r$ and again the inductive hypothesis gives a  directed path from $\ell$ to $k$.  Hence, there is a directed path from $i$ to $k$, as desired.
\item
This equivalence now follows directly from (i) using the definition of $\sco_{G(\mathcal{P}_{\doo})}(k)$ and $\cc{i}$. \qedhere
\end{enumerate}
\end{proof}

As an immediate consequence of the above result, we see that if the set of causes of a variable is non-empty, then at least one of the causes must act as the direct cause.
\begin{corollary}
Let $\mathcal{P}_{\doo}$ be a transitive interventional family where $P_{\doo(i)}$ satisfy the composition property.  For $k\in V$, if $\cause{k}\neq\emptyset$, then $\dcause{k}\neq\emptyset$.
\end{corollary}
\begin{proof}
Let $i\in \cause{k}$. We know that there is at least one directed path from $i$ to $k$. The node adjacent to $k$ on the path is a direct cause by definition.
\end{proof}
In the next example, we illustrate how all the different edges can easily occur.  We remark that we can generate cycles quite easily, but in the standard SCM setting, their existence is non-trivial and requires solvability conditions \cite{bon20}.

\begin{example}[A simple example of an arrow, cycle, and arc]
\label{all-example}
Consider a fixed joint distribution $P$ for random variables $(X_1,X_2,X_3)$, where $X_1$ is not independent of $X_2$, and $X_3$ is independent of $(X_1, X_2)$  Notice we may think of the joint distribution $(X_1, X_2)$ as generated via the following functional equations: $(X_1, X_2=\phi(X_1, U))$ or  $(X_1=\psi(X_2, V), X_2)$, where $\phi$ and $\psi$ are deterministic functions, $U$ and $V$ are uniformly distributed on $[0,1]$.

Thus corresponding to  $\phi$, we have an SCM,  where $X_1  \fra X_2$, and $X_3$ is isolated, and similarly, corresponding to $\psi$, we have an SCM, where $X_2  \fra X_1$, and $X_3$ is isolated. Next, we see how these SCM interact with various interventions, and see the corresponding causal graphs that can be defined.

Via the function $\phi$, we consider standard interventions, where $X_1$, $X_2$, and $X_3$  are replaced with a distributional copies of themselves,  giving the family $(P=P^{\phi}_{\doo(1)}, P^{\phi}_{\doo(2)}, P=P^{\phi}_{\doo(3)} )$  and via the function $\psi$ we consider  standard interventions giving the family $(P^{\psi}_{\doo(1)}, P=P^{\psi}_{\doo(2)}, P=P^{\psi}_{\doo(3)})$;  in both cases, we place an arrow between $1$ and $2$ in the expected direction.  Note that
$$P^{\rm{ind}}:= P^1 \otimes P^2 \otimes P^3= P^{\phi}_{\doo(2)} = P^{\psi}_{\doo(1)}.$$

Consider also the non-standard interventional   family $(P^{\phi}_{\doo(1)}, P^{\psi}_{\doo(2)}, P)$; $1$ is a direct cause of $2$, and $2$ is also a direct cause of $1$, so we obtain a parallel edge.

Finally, to have an arc, we consider the interventional family
$$(P_{\doo(1)}, P_{\doo(2)}, P_{\doo(3)} )= (P^{\mathrm{ind}}, P^{\mathrm{ind}}, P);$$
there are no causes, and no arrows are placed,  but $P_{\doo(3)}$
detects the dependence in $X_1$ and $X_2$, and places an arc between them.  \erk
\end{example}

In the following example, we stress that the generated causal graph heavily depends on the interventional family that is considered, even under standard interventions arising from simple SCM.

\begin{example}[Two different graphs, one underlying distribution]
\label{two-triangles}
Let $\epsilon_1,\epsilon_2,\epsilon_3 $ be independent Bernoulli random variables with parameter $\tfrac{1}{2}$.  Consider the  random variables $X_1 = \epsilon_1$, $X_2 = X_1 + \epsilon_2$, and $X_3 = X_1 + \epsilon_3$; they have a joint distribution $P$, and can be thought of as an SCM  $\mathcal{C}$ with $X_1 \fra X_2$ and $X_1 \fra X_3$.  Thus we can define a intervention family $\mathcal{P}_{\doo}$ corresponding to standard interventions on $\mathcal{C}$; it is not difficult to verify in this case that the resulting casual graph will be the same as the graph for $\mathcal{C}$.

However, it is not difficult to construct another SCM $\mathcal{C}'$ with $X_3' \fra X_2'$, $X_3' \fra X_1'$, and $X_2' \fra X_1'$, where $(X_1', X_2', X_3')$ also have the joint distribution $P$.    Thus we can define another intervention family $\mathcal{P}_{\doo}'$ corresponding to standard interventions on $\mathcal{C}'$; again it is not difficult to verify in this case that the resulting casual graph will be the same as the graph for $\mathcal{C}'$.

We exploited the simple fact that that the joint distribution $P$ does not uniquely determine a corresponding SCM---not even up to adjacency.  Notice that  $\mathcal{C}$ and  $\mathcal{C}'$ do not have the same number of edges.
\erk
\end{example}

The generated graph is indeed a bowless directed mixed graph (BDMG).
\begin{proposition}[The causal graph is BDMG]
\label{prop:parancg}
The causal graph and intervened graphs generated from a transitive interventional family are BDMG.
\end{proposition}
\begin{proof}
The proof is immediate from the definition of intervened causal and causal graphs.
\end{proof}
\begin{remark}
\label{rem:DAG}
If we assume that  $\cc{k}=\{k\}$ for every $k \in V$, so that there is no causal cycle, then the causal graph is a bowless ADMG. If we assume that for every $j,k\in V$ that are not direct causes of each other and every $i$, we have
$$j\ci_{P_{\doo(i)}} k\cd \icause_i(\{j,k\}),$$
then the causal graph does not contain arcs -- in this case it is seen that the interventional family does not detect any  latent variables that cause both $j$ and $k$,  since they are not dependent. Assuming both of these conditions results in a DAG.  \erk
\end{remark}

\subsection{The Markov property with respect to the intervened causal graphs}\label{sec:MarInter}

We can now present a main result of this paper.
\begin{theorem}[Interventional distributions are Markovian to intervened graphs]
\label{thm:markov-int}
Let $\mathcal{P}_{\doo} = \ns{P_{\doo(i)}}_{i \in V}$ be a transitive interventional family.  For each $i \in V$, if  $P_{\doo(i)}$ satisfies the intersection property  and the composition property, then  $P_{\doo(i)}$ is Markovian to the $i$-intervened graph, $G_i(\mathcal{P}_{\doo})$.
\end{theorem}
\begin{proof}
By Theorem \ref{prop:pairMark}, it suffices to show that $P_{\doo(i)}$ satisfies \ref{P} w.r.t.\ $G_i(\mathcal{P}_{\doo})$.
Consider first the case where $i,k\in V$ are non-adjacent in $G_i(\mathcal{P}_{\doo})$.
 We have that $i\notin\pa_{G_i(\mathcal{P}_{\doo})}(k)$; moreover, by definition the definition of  $G_i(\mathcal{P}_{\doo})$, we have that $i\notin\dcause{k}$. Hence,
 Proposition \ref{prop:axcon} implies \ref{P}.

Now we look at the case where $j,k \in V$ are distinct from $i$ and are non-adjacent in $G_i(\mathcal{P}_{\doo})$; in $G_i(\mathcal{P}_{\doo})$,  there are no arrows between $j$ and $k$, and furthermore, from  \eqref{put-arc}, the absence of an arc immediately gives \eqref{P}.
\end{proof}

\subsection{Alternative intervened causal graph}
We defined intervened causal graphs $G_i(\mathcal{P}_{\doo})$ by deciding to place arcs in the graph in places where arrows are missing. We used $j\notci_{P_{\doo(i)}} k\cd \icause_i(\{j,k\})$, which immediately implies parts of the global Markov property as expressed in Theorem \ref{thm:markov-int}. If a latent variable that causes both $j$ and $k$ should always induce dependencies between the two regardless of what to condition on, in order to obtain $G_i(\mathcal{P}_{\doo})$, one can instead place an arc in the causal structure when, for \emph{every} $C\subset V$ such that $i,j,k\notin C$, we have
\begin{equation}
\label{everyC}
j\notci_{P_{\doo(i)}} k\cd C.
\end{equation}
We note that the above definition removes any assumptions related to Markov properties. In order to obtain the same graph, and consequently, Markov property, we will need ordered upward- and downward-stability as the following result shows. Let us denote the intervened causal graph generated in this way by $\tilde{G}_i(\mathcal{P}_{\doo})$.
%
\begin{proposition}
Let $\mathcal{P}_{\doo}$ be an interventional family.  For each $i \in V$, if  $P_{\doo(i)}$ satisfies ordered upward- and downward-stability w.r.t.\ the order induced by $i$-intervened cause, then $\tilde{G}_i(\mathcal{P}_{\doo})=G_i(\mathcal{P}_{\doo})$, where we use \eqref{everyC} in lieu of \eqref{put-arc} for arcs in $\tilde{G}_i(\mathcal{P}_{\doo})$.
\end{proposition}
\begin{proof}
By definition, the arrows are the same in both graphs. We need to show that they have the same arcs. One direction is trivial. We prove the other direction in contrapositive from by assuming $j\ci_{P_{\doo(i)}} k\cd C$ for some $C$. Using  ordered upward-stability, we can add $\icause_i(\{j,k\})$ to the conditioning set, and by applying ordered downward-stability we can remove the other nodes to obtain $j\ci_{P_{\doo(i)}} k\cd \icause_i(\{j,k\})$.
\end{proof}

\section{Observable interventional families}\label{sec:obsint}
In Section \ref{sec:causgr}, the causal structures are completely defined by interventional families; here we provide additional assumptions for the causal graph to be Markovian with respect to an underlying distribution that can be observed.

\subsection{Markov property w.r.t.\ the causal graph}

We say that an interventional family of distributions $\mathcal{P}_{\doo}$ on the state space
$\mathcal{X} = \prod_{i\in V}\mathcal{X}_i$
is   \emph{observable} with respect to an \emph{underlying}  distribution $P$ on $\mathcal{X}$ such that the following axiom holds.
\newpage
\begin{axiom}
\label{eq:ax20}
\hspace{1 cm}
\begin{enumerate}[(a) ]
\item
\label{two-arc}
For  every separable pair $j,k\in V$, in $G(\mathcal{P}_{\doo})$, and every distinct $i\in V$, we have
$$j\ci_{P_{\doo(i)}} k\cd \icause_i(\{j,k\})\Rightarrow j\ci_P k\cd \cause{\{j,k\}}.$$
\item
\label{two-edge}
For  every $k$ and every $i \in \cause{k}$, we have
$$i\ci_{P_{\doo(i)}} k\cd \icause_i(\{i,k\})\Rightarrow i\ci_P k\cd \cause{\{i,k\}}.$$
\end{enumerate}
\end{axiom}

Let us remark if $P$ is a product measure on $\mathcal{X}$, then every interventional family will be observable with respect to $P$; in practice we want to consider underlying distributions that have some relation to the interventional family, such as when $P$ is the distribution of an SCM.

\begin{theorem}[The underlying distribution is Markovian to the causal graph]
\label{thm:markov}
Let $\mathcal{P}_{\doo}$ be a transitive observable interventional family, with underlying joint distribution $P$ that satisfies the intersection property and the composition property.  Then  $P$ is Markovian to the causal graph, $G(\mathcal{P}_{\doo})$.
\end{theorem}
\begin{proof}
By Theorem   \ref{prop:pairMark}, it suffices to verify that $P$ satisfies  pairwise Markov property \eqref{P} w.r.t.\ $G(\mathcal{P}_{\doo})$; moreover, we can assume, without loss of generality that $G(\mathcal{P}_{\doo})$ is maximal since separable pairs do not correspond to any separation statements.    If  $j,k \in V$ are non-adjacent nodes, then  there is no arc between $j$ and $k$ in  $G(\mathcal{P}_{\doo})$.  Hence,  by definition, there exists an $i$ such that $j\ci_{P_{\doo(i)}} k\cd \icause_i(\{j,k\})$, from which Axiom \ref{eq:ax20} \eqref{two-arc} gives the desired result.
\end{proof}
%
%
%
\subsection{Causal graphs for observational distributions}
The arcs in the causal graph $G(\mathcal{P}_{\doo})$ can be generated from the observational distribution $P$ and the causes,  if the  inverse of Axiom \ref{eq:ax20} holds.
Assume $\mathcal{P}_{\doo}$ is an observable interventional family of distributions with respect to an underlying observational measure $P$. We say that $\mathcal{P}_{\doo}$ is a \emph{strongly-observable interventional} family if the following additional axiom holds.

\begin{axiom}
\label{eq:ax2}
\hspace{1 cm}
\begin{enumerate}[(a) ]
\item
\label{three-arc}
For every distinct $i,j,k\in V$, we have
$$j\ci_P k\cd \cause{\{j,k\}}\Rightarrow j\ci_{P_{\doo(i)}} k\cd \icause_i(\{j,k\}).$$
\item
\label{three-edge}
For every $k$ and every $i \in \cause{k}$, we have
$$i\ci_P k\cd \cause{\{i,k\}}   \Rightarrow  i\ci_{P_{\doo(i)}} k\cd \icause_i(\{i,k\}).$$
\end{enumerate}
\end{axiom}

Our next example shows that there are univariate observable interventional families that are not given by standard interventions and that univariate observable interventional families may not satisfy Axiom \ref{eq:ax2}.

\begin{example}[Interventions on joint distributions]
\label{inter-joint}
Suppose $X_1$, $X_2$, and $X_3$ are jointly independent (Bernoulli) random variables, with law $P$, which will serve as an underlying distribution.  Now consider the (non-standard) intervention, where $P_{\doo(1)}$ changes the joint distribution of   $X_2$ and $X_3$ to one of {\em dependence},  but leaves the marginal distributions of $X_2$ and $X_3$, and $X_1$ alone; furthermore, we leave $X_1$ independent of $(X_2,X_3)$.
Although we normally think of  $P_{\doo(1)}$ as an interventional on $X_1$, our general definition allows for somewhat counter-intuitive constructions.

Let $P_{\doo(2)}$ and $P_{\doo(3)}$ be  standard interventions on $X_2$ and $X_3$, respectively, that simply leave the original independent distribution unchanged.   The family $\{P_{\doo(1)}, P_{\doo(2)}, P_{\doo(3)}\}$ satisfies Axioms \ref{prop:trans} and \ref{eq:ax20} trivially, but   Axiom \ref{eq:ax2} is not satisfied. \erk
\end{example}

Define $G(P)$ by adding to the causal structure $S(\mathcal{P}_{\doo})$,
    arcs between $j$ and $k$, i.e., $j\arc k$, for nodes $j$ and $k$ not adjacent by an arrow, if
\begin{equation}
\label{must-latent}
j  \notci_P k  \cd \cause{\{j,k\}}.
\end{equation}
Clearly, \eqref{must-latent} suggests the presence of a latent variable; compare with \eqref{put-arc}.
\begin{proposition}\label{prop:causegrp}
Suppose that $\mathcal{P}_{\doo}$ is a strongly-observable interventional family with the underlying  distribution $P$, and $N\geq 3$, so that there are at least three nodes.   Then $G(P)=G(\mathcal{P}_{\doo})$ when these graphs are maximal.
\end{proposition}
\begin{proof}
By definition, the arrows are the same in both graphs. If there is an arc between $j,k$ in $G(\mathcal{P}_{\doo})$, then by definition of arcs, for all $i\in V$, we have  $j\notci_{P_{\doo(i)}} k\cd \icause_i(\{j,k\})$, hence Axiom \ref{eq:ax2} \eqref{three-arc} implies that $j\notci_P k\cd \cause{\{j,k\}}$, and thus there is an arc between $j,k$ in  $G(\mathcal{P})$.   Conversely, if  $j\notci_P k\cd \cause{\{j,k\}}$,  then Axiom \ref{eq:ax20} \eqref{two-arc} implies the existence of the necessary arc in  $G(\mathcal{P}_{\doo})$ for maximal graphs.
\end{proof}

Note that for strongly-observable interventional families, the existence of the $jk$-arc in $G_i(\mathcal{P}_{\doo})$ and $G(\mathcal{P}_{\doo})$ may differ for non-maximal graphs only when $j,k$ is an inseparable pair. This implies that these two graphs are Markov equivalent.




\subsection{Congruent interventional families}\label{sec:cong}

Consider the interventional families $\mathcal{P}_{\doo}=\{P_{\doo(j)}\}_{j\in V}$ and $\mathcal{Q}_{\doo}=\{Q_{\doo(j)}\}_{j\in V}$ over the same state space $\mathcal{X}$.
It is immediate that if \emph{both} families are strongly observable with respect to a single underlying distribution $P$, and if  the causal graphs $G(\mathcal{P}_{\doo})$ and $G(\mathcal{Q}_{\doo})$ are maximal, then they have the same adjacencies.  Motivated by   Axioms \ref{eq:ax20}  and \ref{eq:ax2}, we say that the families are \emph{congruent} if the have the same causes (see Remark \ref{similar-interv-remark}) and

\begin{enumerate}[1.) ]
\item
\label{two-edge-bb}
For  every $k$ and every $i \in {\rm{cause}}^{\mathcal{P}_{\doo}}(k) = {\rm{cause}}^{\mathcal{Q}_{\doo}}(k)$, we have
$$i\ci_{P_{\doo(i)}} k\cd \icause_i ^{\mathcal{P}_{\doo}} (\{i,k\}) \iff i\ci_{Q_{\doo(i)}} k\cd \icause_i ^{\mathcal{Q}_{\doo}} {\{i,k\}}.$$
\item
\label{three-arc-bb}
For every distinct $i,j,k\in V$, we have
$$j\ci_{P_{\doo(i)}} k\cd \icause_i^{\mathcal{P}_{\doo}}(\{j,k\}) \iff j\ci_{Q_{\doo(i)}} k\cd \icause^{\mathcal{Q}_{\doo}}{\{j,k\}}.$$
\end{enumerate}
We collect our observations in the following proposition.

\begin{proposition}
Consider two interventional families over the same state space.
\begin{enumerate}[1.) ]
\item
The families are congruent if and only if their causal graphs are the same.
\item
If the families  are strongly observable with respect to the same underlying distribution, and  their causal graphs are maximal, then the graphs have the same adjacencies; furthermore, if the families have the same causes, and if the causal graphs are ancestral, then the graphs are the same, and  the families are congruent.
\end{enumerate}
\end{proposition}

\begin{proof}
The first claim is immediate.  For the second claim, it is immediate from  Axioms \ref{eq:ax20}  and \ref{eq:ax2},  that the causal graphs have the same  adjacencies; furthermore, under the ancestral assumption, we can also direct the graph using the causes, which are also the same.  Furthermore, if two nodes are adjacent and not causes of each other, we deduce the presence of an arc.  Hence the graphs are the same, and the families are congruent.
\end{proof}

\subsection{Alternative causal graphs}
As another alternative for causal graphs, one can consider a setting where we relax the criterion for an arc to exist in $G(\mathcal{P}_{\doo})$ by only requiring that a corresponding arc exists in $G_i(\mathcal{P}_{\doo})$ for \emph{some} $i$ (as opposed to \emph{every} $i$).

This graph clearly has more arcs than in the previous setting. 
This new setting leads to something extra to the Markov property above when it comes to conditional independencies related to causal graphs: For a separable pair $(j,k)$, Axiom \ref{eq:ax20} always implies that $j\ci_p k\cd \cause{\{j,k\}}$.  On the other hand, the original setting is more consistent with the idea that the presence of an arc indicates the existence of a latent random variable causing the endpoints.

This also leads to Axiom \ref{eq:ax2} to hold under transitivity, composition, and converse Markov property:
\begin{proposition}
Suppose that $\mathcal{P}_{\doo}$ is a transitive observable interventional family with the underlying  distribution $P$. If $P_{\doo(i)}$ satisfy the composition property, and $P$ satisfies the converse Markov property w.r.t.\ $G(\mathcal{P}_{\doo})$, then Axiom \ref{eq:ax2} holds.
\end{proposition}
\begin{proof}
We need two cases. If $j,k$ (where $j$ may be $i$) are non-adjacent in all $G(\mathcal{P}_{\doo})$ then these are not adjacent in every $G_i(\mathcal{P}_{\doo})$. The right-hand-side of the axiom then holds by Markov property of Theorem \ref{thm:markov-int} under transitivity (which ensures, by Theorem \ref{thm:exch}, that ancestors and causes are the same).

If $j,k$ (where $j$ may be $i$)  are adjacent in $G(\mathcal{P}_{\doo})$, then the left-hand-side never holds under the converse pairwise Markov property.
\end{proof}


\section{Specialization to directed ancestral graphs}\label{sec:anc}
When we assume that the causal graph is ancestral (in fact, directed ancestral), in all the definitions and results $\icause_i(\cdot)$ can be replaced by $\cause{\cdot}$ under intersection and composition. In particular, for the definition of $\dcause{\cdot}$, we have the following:
\begin{proposition}
\label{colo:dcause}
Let $\mathcal{P}_{\doo} = \ns{P_{\doo(i)}}_{i \in V}$ be a transitive interventional family, and suppose that $P_{\doo(i)}$ satisfies the intersection and composition properties for every $i \in V$. If $G(\mathcal{P}_{\doo})$ is ancestral,  then for each  $i\in\cause{k}$, we have that
\begin{equation}
\label{ances-def-dcause}
i\in\dcause{k} \text{ if and only if } i\notci_{P_{\doo(i)}} k\cd \cause{k}\setminus\{i\}.
\end{equation}
\end{proposition}
\begin{proof}
Let $C=\cause{k}\setminus\icause_i(k)$. Since $G(\mathcal{P}_{\doo})$ is ancestral, there are no arcs between members of $C$ and $k$ and no directed paths from $k$ to  members of $C$ in $G(\mathcal{P}_{\doo})$.  We have the separation,
$i\dse_m C \cd \icause_i(k)\setminus\{i\}$ in $G_i(\mathcal{P}_{\doo})$, and   from Markov property given in Theorem \ref{thm:markov-int}, we deduce that
\begin{equation}
\label{above-indep}
i\ci_{P_{\doo(i)}} C   \cd \icause_i(k)\setminus\{i\}.
\end{equation}
\begin{enumerate}[label={}]
\item[($\Leftarrow$) ]
If $i\ci_{P_{\doo(i)}} k\cd \cause{k}\setminus\{i\}$, then by \eqref{above-indep} and contraction property, we obtain $$i\ci_{P_{\doo(i)}} C\cup\{k\}\cd \icause_i(k)\setminus\{i\},$$ which  by decomposition implies  $i \not \in \dcause{k}$.
\item[($\Rightarrow$) ]
If $i \not \in \dcause{k}$, then $i\ci_{P_{\doo(i)}} k\cd \icause_i(k)\setminus\{i\}$; then from  the composition property and \eqref{above-indep}, we obtain $i\ci_{P_{\doo(i)}} C\cup\{k\}\cd \icause_i(k)\setminus\{i\}$, and finally by weak union we have $i\ci_{P_{\doo(i)}} k\cd \cause{k}\setminus\{i\}$. \qedhere
\end{enumerate}
\end{proof}
We note that the implication above is, in fact, only using the global Markov property in one direction (and composition, in addition) in the other.

In addition to this new definition of direct cause, and consequently causal structure, we can define a $jk$-arc in the intervened graph $G(\mathcal{P}_{\doo})$ to exist if
\begin{equation}
\label{put-arc-anc}
j\notci_{P_{\doo(i)}} k\cd \cause{\{j,k\}},
\end{equation}
as the following result shows.
\begin{proposition}
\label{colo:dcausep}
Let $\mathcal{P}_{\doo} = \ns{P_{\doo(i)}}_{i \in V}$ be a transitive interventional family, and suppose that $P_{\doo(i)}$ satisfies the intersection and composition properties for every $i \in V$. If $G(\mathcal{P}_{\doo})$ is ancestral,  then for each  $i,j,k$, we have that
$$j\notci_{P_{\doo(i)}} k\cd \cause{\{j,k\}} \text{ if and only if } j\notci_{P_{\doo(i)}} k\cd \icause_i(\{j,k\}).$$
\end{proposition}
\begin{proof}
The proof is a routine variation of the proof of Proposition \ref{colo:dcause}, where instead of \eqref{above-indep}, we use
\begin{equation*}
k\ci_{P_{\doo(i)}} C   \cd \icause_i(\{j,k\}).
\qedhere
\end{equation*}
\end{proof}
\begin{remark}
    \label{anc-axioms}
Consequently, for ancestral causal graphs, Axioms \ref{eq:ax20} and \ref{eq:ax2}, for strongly observable interventional families, can be written with $\cause{\cdot}$ instead of $\icause_i({\cdot})$. \erk
\end{remark}
\section{Quantifiable interventional families}\label{sec:qauntint}
Although transitive observable interventional families lead to Markov property of the observational distribution to the causal graph, they do not allow us to  measure causal effects from observed data.  Measuring causal effects from observed data is an important consequence of the theory presented here, although it is beyond the scope of this paper. Here we provide a more restricted axiom that allows this possibility.

For technical simplicity, we focus on the case where the graphs are directed ancestral.
We recall from  Section \ref{sec:anc} that under the additional assumptions that the interventional family is transitive with each member  satisfying the intersection and composition properties, all occurrences $\icause(\cdot)$ in the definitions and axioms can be replaced with $\cause{\cdot}$ and $\dcause{\cdot}$ using \eqref{ances-def-dcause}.

\subsection{Axiomatization}
If two measures $P$ and $Q$ have the same null sets, then they are \emph{equivalent}, and we write $P \sim Q$.  For a distribution $P$ with state space $\mathcal{X} = \prod_{i\in V} \mathcal{X}_i$, we say that  the family $\ns{\mathcal{P}_{\doo}, P}$ is \emph{compatible} if for  all distinct $i,k \in V$, we have $P^{\cause{k}\cup \ns{k}}_{\doo(i)}  \sim P^{ \cause{k}\cup \ns{k}}$.

We now say that an interventional family (of distributions) $\mathcal{P}_{\doo}$ is \emph{quantifiable interventional}  if there exists an \emph{underlying}  distribution $P$ with state space $\mathcal{X}$ such that the family $\ns{\mathcal{P}_{\doo}, P}$ is compatible and the following axiom holds.
\begin{axiom}
\label{eq:ax1}
For  all distinct $i,k \in V$, there exist regular conditional probabilities such that their marginals on $\mathcal{X}_k$ satisfy
$$
P^k_{\doo(i)}( \cdot \cd x_{\cause{k} \setminus \ns{i} },x_i)=
\begin{cases}
P^k( \cdot \cd x_{\cause{k}\setminus \ns{i}}),  \text{ if } i \not \in \cause{k} &  \\
  P^k( \cdot \cd x_{\cause{k} \setminus \ns{i} },x_i),  \text{ if }  i  \in \cause{k}
\end{cases}
$$
for all $x_{\cause{k}\setminus \ns{i} } \in \mathcal{X}_{\cause{k} \setminus \ns{i}}$  and all $x_i \in \mathcal{X}_i$.
\end{axiom}
In addition, if  $X=(X_j)_{j\in V}$ is a random vector with distribution $P$, then we say that $X_i$ is a \emph{cause} of $X_k$, if $i\in\cause{k}$, and we write $\cause{X_k}$ to denote the set of the \emph{causes} of $X_k$; and similarly for $\dcause{X_k}$.

When Axiom \ref{eq:ax1} holds, we can drop the occurrence of the $x_i$ in the conditioning set of the  distribution.
\begin{proposition}
\label{prop:dopeq}
Let $\mathcal{P}_{\doo}$ be a quantifiable interventional family with the underlying distribution $P$. For  all distinct $i,k \in V$, there exist regular conditional probabilities such that their marginals on $\mathcal{X}_k$ satisfy
$$P^k_{\doo(i)}( \cdot \cd x_{\cause{k}}) = P^k( \cdot \cd x_{\cause{k}})$$
for all $x_{\cause{k}} \in \mathcal{X}_{\cause{k}}$.
\end{proposition}
\begin{proof}
If $i \in \cause{k}$,
then the result is obvious by  Axiom \ref{eq:ax1}; otherwise, by elementary manipulations, and  Axiom \ref{eq:ax1},
we have
\begin{eqnarray*}
P^k_{\doo(i)}( \cdot \cd x_{\cause{k}}) &=& \int P^k_{\doo(i)}( \cdot \cd x_{\cause{k}}, x_i) dP^i_{\doo(i)}(x_i \cd x_{\cause{k}})  \\
&=&   \int P^k( \cdot \cd x_{\cause{k}}) dP^i_{\doo(i)}(x_i \cd x_{\cause{k}}) \\
&=&   P^k( \cdot \cd x_{\cause{k}}) \int dP^i_{\doo(i)}(x_i \cd x_{\cause{k}})  \\
&=& P^k( \cdot \cd x_{\cause{k}}).
 \qedhere
\end{eqnarray*}
\end{proof}


In the corollary below, we see that Axiom \ref{eq:ax1} imposes strong relationships among different $P_{\doo(i)}$.
\begin{corollary}
\label{prop:pdorel}
Let $\mathcal{P}_{\doo}$ be a quantifiable interventional family. For every distinct $i,j,k\in V$, we have
$$P^k_{\doo(j)}( \cdot \cd x_{\cause{k}})=P^k_{\doo(i)}( \cdot \cd x_{\cause{k}})= P^k( \cdot \cd x_{\cause{k}}).$$
\end{corollary}
\begin{proof}
The first equality is a direct consequence  of Axiom \ref{eq:ax1}. 
 The second equality is Proposition \ref{prop:dopeq}.
\end{proof}
\begin{remark}\label{rem:uinque}
Thus,  given $\mathcal{P}_{\doo}$ and Axiom \ref{eq:ax1}, all $P^k( \cdot \cd x_{\cause{k}})$ are uniquely determined; this will lead to the whole distribution $P$ being uniquely determined in the case where there are no latent variables under certain conditions; see Proposition \ref{coro:unique}.  \erk
\end{remark}
\begin{remark}
By Corollary \ref{prop:pdorel}, Axiom \ref{eq:ax1} imposes strong relationships between different $P_{\doo(i)}$. On the other hand, Axioms \ref{eq:ax20} and \ref{eq:ax2} only deal with conditional independencies implied by interventional and underlying distributions; thus, they clearly do not impose the restrictions that appear in Corollary \ref{prop:pdorel}. This shows that there are many interventional families that are strongly observable but not quantifiable.  \erk
\end{remark}
\begin{lemma}
\label{lem:oldax}
Let $\mathcal{P}_{\doo}$ be a quantifiable interventional family with the underlying distribution $P$ and directed ancestral causal graph $G(\mathcal{P}_{\doo})$.
For $i \in \cause{k}$, we have
\begin{equation}\label{9a}
i\ci_{P_{\doo(i)}} k\cd \cause{k}\setminus\{i\}\Longleftrightarrow i\ci_P k\cd \cause{k}\setminus\{i\}.
\end{equation}
\end{lemma}
\begin{proof}
Observe that the dependence
\begin{equation}
\label{depends-star}
i\notci_{P_{\doo(i)}} k\cd \cause{k}\setminus\{i\}
\end{equation}
is equivalent to the existence of measurable subset $F$ of positive measure under $P_{\doo(i)}^{\cause{k} \setminus \ns{i}}$ such that for every $z \in F$ there exist disjoint  measurable subsets $W^{*}, W^{**} \subset \mathcal{X}_i$ of positive measures under the conditional probability  $P_{\doo(i)}^{k}(\cdot \cd z)$ satisfying the inequality
$$P_{\doo(i)}^{k} (\cdot \cd z,x_i^*)\neq P_{\doo(i)}^{k}(\cdot \cd z,x_i^{**}),$$
for all $ x_i^* \in W^{*}$ and all $x_i^{**} \in W^{**}$; in the case that $\cause{k} \setminus \ns{i}$ is empty, we simply omit reference to the variable $z \in F$.

Axiom \ref{eq:ax1} allows us to interchange $P$ and $P_{\doo(i)}$ and infer that  the dependence \eqref{depends-star} is  equivalent to
$$P_{}^{k} (\cdot \cd z,x_i^*)\neq P_{}^{k}(\cdot \cd z,x_i^{**}),$$
which, with compatibility,  is itself equivalent to dependence
\begin{equation*}
i\notci_P k\cd  \cause{k}\setminus\{i\}.
\qedhere
\end{equation*}
\end{proof}

\begin{proposition}\label{thm:conpair}
Let $\mathcal{P}_{\doo}$ be a transitive quantifiable interventional family with the underlying distribution $P$ and directed ancestral causal graph $G(\mathcal{P}_{\doo})$. Assume that $j$ is a direct cause of $k$. Then
\begin{equation*}
  j\notci_P k\cd \cause{k}\setminus\{j\}. 
\end{equation*}
\end{proposition}
\begin{proof}
Assume $j\in\dcause{k}$. 
Thus, the left-hand side of equivalence  \eqref{9a} in Lemma \ref{lem:oldax} never holds with $i$ being replaced by $j$. Hence, we have the result.
\end{proof}

The following examples may help to illustrate  the need for compatibility in our proof of Lemma \ref{lem:oldax}.

\begin{example}[Supports]
\label{supports-counter}
Consider the random variables, $X_1$ and $X_2$,  where $X_1$ is a random isometry of $\mathbb{R}^2$, chosen uniformly from a set of a rotations, and $X_2 = X_1 \epsilon_2$, where $\epsilon_2$ is a standard bivariate normal random variable that is independent of $X_1$,  and $X_1$ acts on $\epsilon_2$ in the natural way.   It is well-known that $X_1$ is independent of $X_2$.    However, on a standard intervention, if the possibility of translations are included, then $\tilde{X}_1$ would be dependent of $X_2$, and
$1 \in \dcause{2}$.
Notice that  compatibility is not satisfied, since it allows for translations that were a null set under the underlying distribution for $X_1$, which was restricted to rotations.  \erk
\end{example}
With regards to compatibility, it is not enough for two measures to have common support.

\begin{example}[Mutually singular measures with common support]
Let $U$ be uniformly distributed on the unit interval.  Note that the bits  of the binary expansion of $U$ given by sequence $b(U)$ are  iid Bernoulli random variables with parameter $p = \tfrac{1}{2}$.  Working in reverse, we let $V$ be a random variable on the unit interval such that $b(V)$ is a sequence of iid Bernoulli random variables with parameter $p = \tfrac{1}{3}$.    Note that although $U$ and $V$ both have the unit interval as their support, their associated laws are mutually singular.

Consider random variables, $X_1$ and $X_2= \phi(X_1, \epsilon_2)$, where $\epsilon_2$ is uniformly distributed in the unit interval, and independent of $X_1$, and the procedure for generating $X_2$ is as follows.  We examine the  input $X_1$ by taking $b(X_1)$ and then taking the (infinite) sample average; once the sample average $q$  is obtained, we generate a single Bernoulli random variable with parameter $q$.  If the distribution of $X_1$ is fixed to that of either $U$ or $V$ above,  then clearly $X_2$ is independent of $X_1$.  However, with an  intervention on $X_1$ that allows for a mixture of $U$ and $V$, we will have
$1 \in \dcause{2}$. \erk
\end{example}
\subsection{Bivariate-quantifiable interventional families}\label{sec:multint}
As before, let $V=\{1,\ldots,N\}$. Let $\mathcal{P}_{\doo}$ be a quantifiable interventional family. We say that $\mathcal{P}_{\doo}$ is a \emph{bivariate-quantifiable interventional} family (of distributions) 
 if the underlying  distribution $P$ satisfies the following:
\begin{axiom}
\label{eq:ax2n}
For  all distinct $i,j,k \in V$ such that $j\notin\cause{k}$ and $k\notin\cause{j}$, there exist regular conditional probabilities such that their marginals on $\mathcal{X}_A$ satisfy
$$
P^{\{j,k\}}_{\doo(i)}( \cdot \cd x_{\cause{\{j,k\}} \setminus \ns{i} },x_i)=
\begin{cases}
P^{\{j,k\}}( \cdot \cd x_{\cause{\{j,k\}}}),  \text{ if } i \not \in \cause{\{j,k\}} &  \\
  P^{\{j,k\}}( \cdot \cd x_{\cause{\{j,k\}} \setminus \ns{i} },x_i),  \text{ if }  i  \in \cause{\{j,k\}}
\end{cases}
$$
for all $x_{\cause{\{j,k\}}\setminus \ns{i} } \in \mathcal{X}_{\cause{\{j,k\}} \setminus \ns{i}}$  and all $x_i \in \mathcal{X}_i$.
\end{axiom}

We note that a bivariate version above does not hold in general, and in particular, may not hold when the intervention $i$ takes place in \emph{between} $j$ and $k$.
\begin{example}
Consider again the  simple SCM, where $1 \fra 2 \fra 3$ and $1 \fra 3$.   In a discrete setting, we have
$$P^{1,3}_{\doo(2)}(x_1, x_3 \cd x_2) = \mathbb{P}(X_1 = x_1, \phi_3(X_1, x_2, \epsilon_3) = x_3),$$
but
$$P_{\mathcal{C}}^{1,3}(x_1, x_3  \cd x_2) =  \mathbb{P}(X_1 = x_1, \phi_3(X_1, x_2, \epsilon_3) = x_3 | X_2 = x_2),$$
and in general, the expressions will not be equal.
\erk
\end{example}

%

We have a similar result to that of Proposition \ref{prop:dopeq} for bivariate observable interventional families.

\begin{proposition}
\label{prop:dopeqmult}
Let $\mathcal{P}_{\doo}$ be a bivariate-quantifiable interventional family with the underlying distribution $P$. Then
for  all distinct $i,j,k \in V$ such that $j\notin\cause{k}$ and $k\notin\cause{j}$,  there exist  regular conditional probabilities such that their marginals on $\mathcal{X}_{\ns{j,k}}$ satisfy
$$P^{\{j,k\}}_{\doo(i)}( \cdot \cd x_{\cause{\{j,k\}}}) = P^{\{j,k\}}( \cdot \cd x_{\cause{\{j,k\}}})$$
for all $x_{\cause{\{j,k\}}} \in \mathcal{X}_{\cause{\{j,k\}}}$.
\end{proposition}

\begin{proof}
The proof is a simple routine modification of the proof of   Proposition \ref{prop:dopeq} for the univariate case.
\end{proof}

For bivariate-quantifiable interventional families, we have the following important conditional independence result.
\begin{corollary}[Conditional independence of a disjoint pair on intervention]
\label{prop:intindstrict}
Let $\mathcal{P}_{\doo}$ be a bivariate-quantifiable interventional family. Then, for distinct $i,j,k\in V$ such that $j\notin\cause{k}$ and $k\notin\cause{j}$, we have that
$$j\ci_{P_{\doo(i)}} k\cd \cause{\{j,k\}} \Longleftrightarrow j\ci_P k\cd \cause{\{j,k\}}.$$
\end{corollary}
\begin{proof}
The result is an immediate consequence of Proposition \ref{prop:dopeqmult}.
\end{proof}
We now provide conditions under which quantifiable families are observable.
\begin{theorem}[Bivariate-quantifiable interventional families are strongly-observable interventional]
\label{thm:biquanttoobint}
Let $\mathcal{P}_{\doo}$ be a transitive bivariate-quantifiable interventional family with the underlying distribution $P$ and a directed ancestral causal graph $G(\mathcal{P}_{\doo})$.  If $P_{\doo(i)}$ satisfies the intersection and composition properties for every $i \in V$,  then  $\mathcal{P}_{\doo}$ is strongly observable.
\end{theorem}
\begin{proof}
Note that by Remark \ref{anc-axioms}, we may replace $\icause({\cdot})$ with $\cause{\cdot}$ in the axioms.

Axioms \ref{eq:ax20} \eqref{two-edge} and \ref{eq:ax2} \eqref{three-edge} now follow from Lemma \ref{lem:oldax} and transitivity.

We prove Axioms \ref{eq:ax20} \eqref{two-arc} and \ref{eq:ax2} \eqref{three-arc}. If $j\notin\cause{k}$ and $k\notin\cause{j}$,  then the result follows directly from Corollary \ref{prop:intindstrict}. Hence, assume, without loss of generality, that $j\in\cause{k}$. In this case, $\cause{\{j,k\}}=\cause{k}\setminus\{j\}$.   By conditioning, for measurable $J \subseteq \mathcal{X}_j$ and $K \subseteq \mathcal{X}_k$, we have
$$P^{\{j,k\}}_{\doo(i)}(J,K \cd x_{\cause{k}\setminus \{j\}}) = \int_K\int_J dP^{k}_{\doo(i)}( x_k \cd x_{\cause{k}\setminus \{j\}}, x_j)  dP^{j}_{\doo(i)}( x_j \cd x_{\cause{k}\setminus \{j\}})$$
and
$$P^{\{j,k\}}(J,K \cd x_{\cause{k}\setminus \{j\}}) = \int_K\int_J dP^{k}( x_k \cd x_{\cause{k}\setminus \{j\}}, x_j)  dP^{j}( x_j \cd x_{\cause{k}\setminus \{j\}})$$
By Proposition \ref{prop:dopeq}, we have
$$P^{k}_{\doo(i)}( \cdot \cd x_{\cause{k}\setminus \{j\}}, x_j) = P^{k}( \cdot \cd x_{\cause{k}\setminus \{j\}}, x_j).$$
Hence, together with  conditional independence $j\ci_{P_{\doo(i)}} k\cd \cause{k} \setminus \ns{j}$, we have the disintegration
\begin{eqnarray*}
P^{\{j,k\}}(J,K \cd x_{\cause{k}\setminus \{j\}}) &=& \int_K\int_J dP^{\{k\}}_{\doo(i)}( x_k \cd x_{\cause{k}\setminus \{j\}}, x_j)  dP^{\{j\}}( x_j \cd x_{\cause{k}\setminus \{j\}}) \\
&=&
\int_K dP^{k}_{\doo(i)}( x_k \cd x_{\cause{k}\setminus \{j\}}) \int_J   dP^{j}( x_j \cd x_{\cause{k}\setminus \{j\}}) \\
&=& P^{k}_{\doo(i)}(K \cd x_{\cause{k}\setminus \{j\}}) P^{j}( J \cd x_{\cause{k}\setminus \{j\}}),
\end{eqnarray*}
 and setting $J=\mathcal{X}_j$ we have
$$P^{k}(K \cd x_{\cause{k}\setminus \{j\}}) = P^{k}_{\doo(i)}(K \cd x_{\cause{k}\setminus \{j\}})$$
from which we have the   conditional independence $j\ci_{P} k\cd \cause{k} \setminus \ns{j}$; the other required direction of the proof is similar.
\end{proof}
As a direct consequence  of Theorem \ref{thm:markov}, we have the following Markov property.
\begin{corollary}
\label{use-later-in-dag}
Let $\mathcal{P}_{\doo}$ be a transitive bivariate-quantifiable interventional family, with underlying joint distribution $P$ and directed ancestral causal graph $G(\mathcal{P}_{\doo})$. Assume that $P$ and $P_{\doo(i)}$, for every $i \in V$, satisfy the intersection property and the composition property. It then holds that $P$ is Markovian to the causal graph, $G(\mathcal{P}_{\doo})$.
\end{corollary}
\begin{corollary}\label{coro:conpair}
Let $\mathcal{P}_{\doo}$ be a transitive quantifiable interventional family with the underlying distribution $P$ and directed ancestral causal graph $G(\mathcal{P}_{\doo})$.   If $P_{\doo(i)}$ satisfies the intersection and composition properties for every $i \in V$, then $P$ satisfies the converse pairwise Markov property w.r.t.\ $G(\mathcal{P}_{\doo})$.
\end{corollary}
\begin{proof}
The case where there is an arrow from $j$ to $k$ is Proposition \ref{thm:conpair}.

The required dependence,  $j\notci_{P_{\doo(i)}} k\cd \cause{\{j,k\}}$,  where there is an arc between $j$ and $k$ in $G(\mathcal{P}_{\doo})$ follows from Theorem \ref{thm:biquanttoobint} and Proposition \ref{prop:causegrp}.
\end{proof}
\subsection{Uniqueness of observable distribution}
We can now prove the uniqueness of $P$ under Axiom \ref{eq:ax1} for DAGs.
\begin{proposition}
\label{coro:unique}
Let $\mathcal{P}_{\doo}$ be a transitive bivariate-quantifiable interventional family, with underlying joint distribution $P$. Assume that $P$ and $P_{\doo(i)}$, for every $i \in V$, satisfy the intersection property and the composition property.  If $G(\mathcal{P}_{\doo})$ is a DAG then $P$ is unique.
\end{proposition}
\begin{proof}
Notice that the Markov property in Corollary  \ref{use-later-in-dag} leads to a factorization of $P$ for DAGs.
For measurable $F \subset \mathcal{X} = \prod_{i \in V} \mathcal{X}_i$,
we have
\begin{eqnarray*}
P(F) &=& \int_A dP(x) \\
&=&\int_F\prod_{k\in V}dP^k(x_k \cd x_{>k}) \\
&=& \int_F\prod_{k\in V}dP^k(x_k\cd x_{\cause{k}}),
\end{eqnarray*}
where $>k$ is the set of all nodes larger than $k$, which is obtained from  a valid ordering of nodes (see Corollary \ref{coro:order}), and the last equality uses the independence $k\ci_P (>k)\setminus{\cause{k}}\cd \cause{k}$, implied by Markov property.   The desired uniqueness now follows from Remark \ref{rem:uinque}.
\end{proof}
Notice that, in the case where there is an arc, the uniqueness does not hold. One can simply consider an $ij$-arc to observe that, using Axiom \ref{eq:ax1}, only the marginal distributions of $X_i$ and $X_j$ are determined.

\section{Specialization to structural causal models}\label{sec:scm}
In this section, we relate the standard intervention on SCMs to the setting presented in this paper.

Let $\mathcal{C}$ be an SCM with random vector $X$ taking values on $\mathcal{X} = \prod_{i \in V} \mathcal{X}_i$ with  joint distribution $P_{\mathcal{C}}$, and associated graph $G_{\mathcal{C}}$. Consider again the standard intervention in SCMs, where intervention on $i \in V$ replaces the equation $X_i=\phi_i(X_{\pa_G(i)},\epsilon_i)$ by $X_i=\tilde{X}_i$, where $\tilde{X}_i$ is independent of all other noises. In the setting of this paper, the new system of equations after intervening on $i$ yields the joint distribution
$P_\mathcal{C}{\scriptstyle[{\doo(i)=\tilde{X}_i}]}$, and consequently one obtains the family of distributions $\mathcal{P}_{\mathcal{C}}{\scriptstyle[{\doo=\tilde{X}}]} :=\{P_{\mathcal{C}}{\scriptstyle[{\doo=\tilde{X}_1}]}, \ldots, P_{\mathcal{C}}{\scriptstyle[{\doo=\tilde{X}_N}]}\}.$
\begin{remark}
The definition of the set $\cause{i}$, defined for $\mathcal{P}_{\mathcal{C}}{\scriptstyle[{\doo=\tilde{X}}]}$, in the setting of this paper is identical to the definition of the set of ``causes'' of $i$ in the SCM setting \cite{pea09,pet17}. \erk
\end{remark}
\begin{remark}
\label{peters-thing}
Note that under some weak, but technical assumptions (in the sense of compatibility), we suspect it is  possible to show that  $\cause{A}$ is invariant under the choice of $\tilde{X}$; see also Proposition 6.13 in \cite{pet17},  which has technical counterexamples.

Therefore, under invariance, if we are only interested in the causal structure, we can simply refer to some canonical  family, where intervention $\tilde{X}_i$ has the same distribution as $X_i$, which we denote by
$\mathcal{P}_{\doo}(\mathcal{C}) = \{P_{\doo(1)}, \ldots, P_{\doo(N)}\}$. \erk
\end{remark}
\subsection{SCMs and interventional families}
The first question that needs to be addressed is when $\mathcal{P}_{\doo}(\mathcal{C})$ satisfies different axioms and key assumptions related to interventional families. It is well known that cancellations may occur in SCMs so that cause is not transitive (as required in Axiom \ref{prop:trans}). We do not provide conditions for SCMs such that cause is transitive-- the main results in this section do not require transitivity of cause.
The following example shows that standard interventions on SCMs do not lead to quantifiable interventional families either.
\begin{example}\label{ex:doubledo}
Consider the collider $X_1 \fra X_3 \fla X_2$, where $X_3 = X_1 \oplus X_2 \bmod 2$.  Consider the underlying joint distribution $P_{\mathcal{C}}$ where $X_1$ is Bernoulli with parameter $p_1=\tfrac{1}{100}$ and $X_2$ is Bernoulli with parameter $p_2=\tfrac{1}{2}$.   Consider the standard interventions where $p_1\to \tfrac{1}{2}$ and $p_2 \to \tfrac{1}{100}$.  Although these are standard interventions, the resulting family does not satisfy Axiom \ref{eq:ax1}.  Observe that  $X_2$ is a direct cause of $X_3$, but $X_1$ is not a cause of $X_3$. However, $P_{\doo(1)}(x_3=1| x_2=1, x_1=1) = 0$
and
$P(x_3=1 | x_2=1) = P(x_1=0) = \tfrac{99}{100}.$  We also see that the conditions of Theorem \ref{prop:intind} are not satisfied.
\erk
\end{example}
We will need to introduce a concept related to faithfulness on the edge level.   We say that $\mathcal{P}_{\doo}(\mathcal{C})$ satisfies the \emph{edge-cause} condition w.r.t.\ $G_{\mathcal{C}}$  if an arrow from $i$ to $j$ in $G_{\mathcal{C}}$ implies that $i\in\dcause{j}$, i.e., $i\notci_{P_{\doo(i)}}j$ and $i\notci_{P_{\doo(i)}}j\cd \icause_i(j)\setminus\{i\}$.

In Section \ref{section-jlmr},   we  will discuss simple conditions on the SCM that imply the edge-cause condition. In particular, it is easy to see that if $P_{\doo(i)}$  are faithful to the intervened graphs $(G_{\mathcal{C}})_i$, where the upcoming arrows and arcs to $i$ are removed, the edge-cause condition is satisfied. The same can be said with the weaker condition of adjacency-faithfulness of $P_{\doo(i)}$ and $(G_{\mathcal{C}})_i$.


\begin{proposition}[Ancestors and causes]
\label{prop:causeanscm}
For an  SCM  $\mathcal{C}$ and the family $\mathcal{P}_{\doo}(\mathcal{C})$, we have that
\begin{equation}
\label{cause-sub}
\dcause{k} \subseteq \pa_{G_{\mathcal{C}}}(k) \text { and }
\cause{k} \subseteq \an_{G_{\mathcal{C}}}(k),
\end{equation}
for every $k\in V$. In addition, if $\mathcal{P}_{\doo}(\mathcal{C})$ satisfies the edge-cause condition w.r.t.\ $G_{\mathcal{C}}$ then
\begin{equation}
\label{cause-par}
\dcause{k}=\pa_{G_{\mathcal{C}}}(k) \text { and }
\pa_{G_{\mathcal{C}}}(k) \subseteq \cause{k}.
\end{equation}
Moreover, $\mathcal{P}_{\doo}(\mathcal{C})$ is transitive if and only if
\begin{eqnarray}
\label{cause-an}
\cause{k}&=&\an_{G_{\mathcal{C}}}(k).
\end{eqnarray}
\end{proposition}
\begin{proof}
By the structure of the SCM,  if $j\notin\an_{G_{\mathcal{C}}}(k)$,  then an intervention $\tilde{X}_j$ on  $X_j$ does not affect the distribution of $X_k$, so that $\tilde{X}_j$ would be independent of $X_k$, and would not be a cause of $k$.    Similarly, if  $j\notin\pa_{G_{\mathcal{C}}}(k)$, then  given $\icause_j{k}\setminus\{j\}$ an intervention $\tilde{X}_j$ on  $X_j$ does not affect the distribution of $X_k$, and thus $j$ cannot be a direct cause of $k$.  Hence we have established \eqref{cause-sub}, and \eqref{cause-par} follows from the definition of the edge-cause condition.

To prove \eqref{cause-an}, we note that we already have $\dcause{k}=\pa_{G_{\mathcal{C}}}(k)$. Transitivity of the interventional family and $\cause{k}\subseteq\an_{G_{\mathcal{C}}}(k)$ imply that $\cause{k}=\an_{G_{\mathcal{C}}}(k)$; the other direction follows from the fact that the set of ancestors is transitive.
\end{proof}

 The following example shows that the inequality may be strict in \eqref{cause-sub}; and  \eqref{cause-par} does not hold without the edge-cause condition.
\begin{example}[Independence and cause in an SCM]\label{ex:col2}
Consider the collider $$X_1 \fra X_3 \fla X_2,$$ where $X_3 = X_1 \oplus X_2 \bmod 2$.  Consider the underlying joint distribution $P_{\mathcal{C}}$ where $X_1$ is Bernoulli with parameter $p_1=\tfrac{1}{2}$ and $X_2$ is Bernoulli with parameter $p_2=\tfrac{1}{2}$.   Consider the standard interventions where nothing happens: $p_1\to \tfrac{1}{2}$ and $p_2 \to \tfrac{1}{2}$.   Then $X_3$ has no causes, and thus its set of ancestors is not equal to its causes.  However, it is easy to verify that Axiom \ref{eq:ax1} is satisfied.   \erk
\end{example}

We also recall that the joint distribution $P_{\mathcal{C}}$ of a structural causal model $\mathcal{C}$ is Markovian to $G_{\mathcal{C}}$. We need a corresponding result to the Markov property of the joint distribution of the SCMs for the intervened distribution.
\begin{lemma}\label{lem:scmintMark}
Let $\mathcal{C}$ be an SCM.  For each $i\in V$, its intervened distribution $P_{\doo(i)}$ is Markovian to $(G_{\mathcal{C}})_i$ and $G_{\mathcal{C}}$.
\end{lemma}
\begin{proof}
Notice that intervention on $i$ yields another SCM with joint distribution $P_{\doo(i)}$ and intervened graph $G_i$, where $G_i$ is obtained from $G_{\mathcal{C}}$ by removing all the parents of $i$ and all the arcs with one endpoint $i$. Therefore, $P_{\doo(i)}$ is Markovian to $G_i$. Since, $G_i$ is a strict subgraph of $G_{\mathcal{C}}$ (with the same node set), so that  a separation in the original graph is a separation in the intervened graph,  we conclude that $P_{\doo(i)}$ is Markovian to $G_{\mathcal{C}}$.
\end{proof}

\begin{theorem}[Strongly observable SCMs]
\label{coro:scmstrobs}
Let $\mathcal{C}$ be an SCM associated to graph $G_{\mathcal{C}}$. 
 Assume that $P_{\mathcal{C}}$ satisfies the converse pairwise Markov property w.r.t.\ $G_{\mathcal{C}}$.  Also assume that $\mathcal{P}_{\doo}(\mathcal{C})$ satisfies the edge-cause condition w.r.t.\ $G_{\mathcal{C}}$.   Then  $\mathcal{P}_{\doo}(\mathcal{C})$ is a strongly observable interventional family if $G_{\mathcal{C}}$ is ancestral. In addition, if $\mathcal{P}_{\doo}(\mathcal{C})$ is transitive,
 then the result holds for BDMGs.
\end{theorem}
\begin{proof}
To prove Axiom \ref{eq:ax20} \eqref{two-arc} in the ancestral case, by Lemma \ref{lem:squeeze}, for a separable pair $(j,k)$,  we have that $j \dse_{m} k\cd A$, for any subset $A$ of the ancestors of $j$ and $k$ containing their parents.  Under the edge-cause condition, $\cause{\{j,k\}}$ contains the parents. Hence, by the Markov property of SCMs, we have that $j\ci_P k\cd \cause{\{j,k\}}$.

To prove Axiom \ref{eq:ax20} \eqref{two-edge}, in the ancestral case,  suppose that $i \in \cause{k}$. By Proposition \ref{prop:causeanscm}, we have that $i \in \an_{G_{\mathcal{C}}}(k)$. By the edge-cause condition, we may assume that there is no arrow from $i$ to $k$. There is also no arc between $i$ and $k$, since causes are ancestors. Since in ancestral graphs, inseparable pairs of nodes cannot be ancestors of each other, we conclude that $i,k$ are separable.   Again by using Lemma \ref{lem:squeeze}  and the same argument as for  Axiom \ref{eq:ax20} \eqref{two-arc},   the result follows from the Markov property for SCMs.

If we further assume that $\mathcal{P}_{\doo}(\mathcal{C})$ is transitive then, by Proposition \ref{prop:causeanscm}, we have that $\cause{\{i,k\}}=\an(\{i,k\})$, and the results follow without the assumption of the graph being ancestral.

To prove Axioms \ref{eq:ax2} \eqref{three-arc} and \eqref{three-edge} , we will consider two cases. If $j,k \in V$ are not adjacent in  $G_{\mathcal{C}}$, then the result follows from Markov property of $P_{\doo(i)}$ for SCMs, given in  Lemma \ref{lem:scmintMark}. If $j$ and $k$ are adjacent, then the result follows from the converse pairwise Markov property.
\end{proof}

\subsection{Causal graphs and graphs associated to SCMs}
\label{section-proof}
In this subsection, we present the ultimate relationship between interventions in the SCM setting and the setting in this paper, which relates the ``true causal graph'' $G_{\mathcal{C}}$ with the causal graph $G(\mathcal{P}_{\doo})$ defined in this paper. We first need some lemmas.
\begin{lemma}\label{lem:causean}
Let $\mathcal{C}$ be an SCM.  We have that
$\pa_{G(\mathcal{P}_{\doo}(\mathcal{C}))}(k) \subseteq \pa_{G_{\mathcal{C}}}(k),$
for every $k\in V$. 
 In addition, if $\mathcal{P}_{\doo}(\mathcal{C})$ satisfies the edge-cause condition w.r.t.\ its associated graph $G_{\mathcal{C}}$,  then
 $\pa_{G(\mathcal{P}_{\doo}(\mathcal{C}))}(k)=\pa_{G_{\mathcal{C}}}(k).$
\end{lemma}
\begin{proof}
Let $j\in \pa_{G(\mathcal{P}_{\doo}(\mathcal{C}))}(k)$. This arrow is a direct cause by construction; hence, by Proposition \ref{prop:causeanscm}, it follows that  $j\in \pa_{G_{\mathcal{C}}}(k)$.

To prove the second statement, if $i \in \pa_{G_{\mathcal{C}}}(k)$, then by the edge-cause condition, we have $i\in\dcause{j}$. Hence, by construction,
 $i\in \pa_{ G(\mathcal{P}_{\doo}(\mathcal{C}))}(j)$.
\end{proof}
To prove a part of the next main result that deals with directed ancestral graphs, we need the following lemma.
\begin{lemma}\label{lem:ancausescm}
Let $\mathcal{C}$ be an SCM, and assume that  $\mathcal{P}_{\doo}(\mathcal{C})$ satisfies the edge-cause condition w.r.t.\ its associated directed ancestral graph $G_{\mathcal{C}}$. Then for every $i,j \in V$, we have
$$i\ci_{P_{\mathcal{C}}}j\cd \cause{\{i,j\}}\Rightarrow i\ci_{P_{\mathcal{C}}}j\cd \an_{\mathcal{C}}(\{i,j\}).$$
\end{lemma}
\begin{proof}
Assume that $i\ci_{P_{\mathcal{C}}}j\cd \cause{\{i,j\}}$; in addition, without loss of generality, choose $j$ not to be a descendent of $i$.
Let $K:=\an_{G_{\mathcal{C}}}(\{i,j\})\setminus\cause{\{i,j\}}$, and $k\in K$. By edge-cause condition,  $k\notin\pa(i)\cup\pa(j)$, and since $G_{\mathcal{C}}$ is ancestral, we know that there is no arc between $k$ and $\{i,j\}$.    Hence we have the separation,  $i\dse_m k\cd \cause{\{i,j\}}\cup\{j\}$; moreover,  since $P_{\mathcal{C}}$ is Markovian to $G_{\mathcal{C}}$, we have the independence,  $i\ci_{P_{\mathcal{C}}} k\cd \cause{\{i,j\}}\cup\{j\}$.

The contraction property implies  that $i\ci_{P_{\mathcal{C}}} \{j,k\}\cd \cause{\{i,j\}}$, and  an inductive argument with similar reasoning, we obtain that
$$i\ci_{P_{\mathcal{C}}} [\{j\}\cup K]  \cd  \cause{\{i,j\}}.$$
Finally,  by weak union, we have $i\ci_{P_{\mathcal{C}}}j\cd \an(\{i,j\}).$
%
\end{proof}
\begin{theorem}[Equality of causal and SCM graphs]
\label{main-SCM-comp}
Consider a structural causal model $\mathcal{C}$ with the joint distribution $P_{\mathcal{C}}$. If $\mathcal{P}_{\doo}(\mathcal{C})$ satisfies the edge-cause condition, and $P_{\mathcal{C}}$ satisfies the converse pairwise Markov property w.r.t.\ the maximal directed ancestral graph $G_{\mathcal{C}}$,  then  $G(\mathcal{P}_{\doo}(\mathcal{C}))=G_{\mathcal{C}}$. In addition, if  $\mathcal{P}_{\doo}(\mathcal{C})$ is transitive, then the result holds for maximal BDMGs. Moreover, without the maximality assumption,
it holds that
$G(\mathcal{P}_{\doo}(\mathcal{C}))$ and $G_{\mathcal{C}}$ are Markov equivalent.
\end{theorem}
\begin{proof} 
By Lemma \ref{lem:causean}, we have that $\pa_{G(\mathcal{P}_{\doo}(\mathcal{C}))}(k)=\pa_{G_{\mathcal{C}}}(k)$. Therefore, in order to show $G(\mathcal{P}_{\doo}(\mathcal{C}))=G_{\mathcal{C}}$, it is enough to prove that $G(\mathcal{P}_{\doo}(\mathcal{C}))$ and $G_{\mathcal{C}}$ have the same arcs. We will check pairs of nodes with no edge between them.

Since $\mathcal{P}_{\doo}(\mathcal{C})$ is strongly observable (Theorem \ref{coro:scmstrobs}), by Proposition \ref{prop:causegrp}, there is no edge between $i$ and $j$ in $G(\mathcal{P}_{\doo}(\mathcal{C}))$ if and only if $i\ci_{P_{\mathcal{C}}}j\cd \cause{\{i,j\}}$;  assume that this conditional independence statement holds. If $G_{\mathcal{C}}$ is ancestral then by Lemma \ref{lem:ancausescm}, we have that $i\ci_{P_{\mathcal{C}}}j\cd \an(\{i,j\})$.

If $G_{\mathcal{C}}$ is not ancestral, then using transitivity, and Proposition \ref{prop:causeanscm},  $\an(\cdot)$ is replaced with $\cause{\cdot}$. Now. in both cases, the converse pairwise Markov property implies that there is no edge between $i$ and $j$ in $G_{\mathcal{C}}$.

Now consider the other direction: Let $i \not \sim j$ in $G_{\mathcal{C}}$. If $G_{\mathcal{C}}$ is ancestral, then, since it is also maximal, and $\cause{\{i,j\}}\subseteq\an_{G_{\mathcal{C}}}{\{i,j\}}$, by Lemma \ref{lem:squeeze}, we have $i\dse_m j\cd \cause{\{i,j\}}$. Then the global Markov property for SCMs implies $i\ci_{P_{\mathcal{C}}}j\cd \cause{\{i,j\}}$.

If $G_{\mathcal{C}}$ is not ancestral, then by Proposition \ref{prop:Markpair}, with $\an(\cdot)$ replaced with $\cause{\cdot}$ (using Proposition \ref{prop:causeanscm}), it holds that $i\ci_{P_{\mathcal{C}}}j\cd \cause{\{i,j\}}$. This implies that there is no arc between $i$ and $j$ in $G(\mathcal{P}_{\doo}(\mathcal{C}))$.
%

Finally, we note that if the graph is not maximal, the only pairs of nodes $(i,j)$ for which $i\ci_{P_{\mathcal{C}}}j\cd \cause{\{i,j\}}$ may not hold are inseparable pairs. However, the lack of an edge between inseparable pairs in  $G_{\mathcal{C}}$, and simultaneously the existence of it in $G(\mathcal{P}_{\doo}(\mathcal{C}))$, do not violate Markov equivalence.
\end{proof}
The following example shows that we, indeed, require to assume that $G_{\mathcal{C}}$ is maximal.
\begin{example}
Consider the non-maximal graph of Figure \ref{fig:nonmax2} to be $G_{\mathcal{C}}$. Assume standard intervention and also faithfulness of $P_{\mathcal{C}}$ and $G_{\mathcal{C}}$.  In $G(\mathcal{P}_{\doo}(\mathcal{C}))$, there exists an arc between $j$ and $k$ since no matter what one intervenes on, $j$ and $k$ always stay dependent given any conditioning set. Notice that here we require two discriminating paths between $j$ and $k$. If there were only one discriminating path between $j$ and $k$ (for example with no $h'$ and $\ell'$) then by intervening on any node on the discriminating path (e.g., $h$), one obtains the required independence $j\ci_{P_{\doo(h)}} k\cd \icause_h(\{j,k\})$. \erk
\begin{figure}[htb]
\centering
\begin{tikzpicture}[node distance = 8mm and 8mm, minimum width = 6mm]
    \begin{scope}
      \tikzstyle{every node} = [shape = circle,
      font = \scriptsize,
      minimum height = 6mm,
      inner sep = 2pt,
      draw = black,
      fill = white,
      anchor = center],
      text centered]
      \node(l) at (0,0) {$\ell$};
      \node(h) [right = of l] {$h$};
      \node(k) [below right = of h] {$k$};
      \node(j) [below left = of l] {$j$};
      \node(l') [below right = of j] {$\ell'$};
      \node(h') [below left = of k] {$h'$};
    \end{scope}
		
    \begin{scope}[->, > = latex']
    \draw (l) to [bend right]  (k);
    \draw (h) to [bend left]  (j);
    \draw (l') to [bend left]  (k);
    \draw (h') to [bend right]  (j);
    \end{scope}
    \begin{scope}[<->, > = latex']
    \draw (j) -- (l);
    \draw (l) -- (h);
    \draw (h) -- (k);
    \draw (j) -- (l');
    \draw (l') -- (h');
    \draw (h') -- (k);
    \end{scope}

    \end{tikzpicture}
		\caption{{\footnotesize A (non-maximal) graph associated to an SCM}}\label{fig:nonmax2}
		\end{figure}
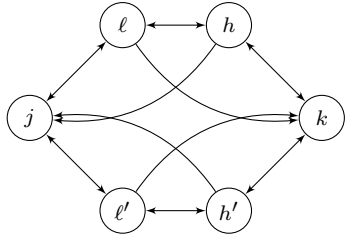
\end{example}


\begin{corollary}
Let $\mathcal{C}$ be an SCM. If its joint distribution $P_{\mathcal{C}}$ and its associated bowless directed mixed graph $G_{\mathcal{C}}$ are faithful, and its intervened distribution $P_{\doo(i)}$ and its associated intervened graphs $(G_{\mathcal{C}})_i$ are faithful, for every $i$, then $G(\mathcal{P}_{\doo}(\mathcal{C}))=G_{\mathcal{C}}$.
\end{corollary}
\begin{proof}
The result follows from the fact that faithfulness of intervened distributions and graphs implies the edge-cause condition and faithfulness of the underlying distribution and graph implies the converse pairwise Markov property. 
\end{proof}

\subsection{SCMs and quantifiable interventional families}
In this subsection, we assume that the graph $G_{\mathcal{C}}$ associated to an SCM is a directed ancestral graph in order to follow the theory presented in Section \ref{sec:qauntint}.

Although  Axiom \ref{eq:ax1} may not in general be satisfied, we do have the following observations. 

\begin{lemma}
\label{scm-fact}
Let $\mathcal{C}$ be an SCM associated to a directed ancestral graph $G_{\mathcal{C}}$.  Let $i \in V$,  $k \in V \setminus \ns{i}$, and  $A \subseteq V$ be such that $\pa(k) \subseteq A \subseteq \an_{G_{\mathcal{C}}}(k)$.
There exist regular conditional probabilities such that their marginals on $\mathcal{X}_k$ satisfy
\begin{equation*}
 P^k(\cdot  \cd x_A) =\mathbb{P}(\phi_k(x_ {\pa_{G_{\mathcal{C}}}(k)} , \epsilon_k) \in \cdot) \text{ \ for all  $x_A \in \mathcal{X}_A$},
\end{equation*}
\begin{equation*}
P^k_{\doo(i)}( \cdot \cd x_{A}) = P^k(\cdot  \cd x_A) \text{ \ for all  $i  \in \an_{G_{\mathcal{C}}}(k)$,    $x_A \in \mathcal{X}_A$},
\end{equation*}
and
\begin{equation*}
P^k_{\doo(i)}( \cdot \cd x_{A}, x_i) = P^k(\cdot  \cd x_A) \text{ \ for all  $i \not \in \an_{G_{\mathcal{C}}}(k)$, $x_i \in \mathcal{X}_i$,    $x_A \in \mathcal{X}_A$}.
\end{equation*}
\end{lemma}

\begin{proof}
 Let $K \subseteq \mathcal{X}_k$ be measurable.   Then for $x_{A}  \in \mathcal{X}_{A}$, we have
\begin{eqnarray*}
\mathbb{P}(X_k \in K| X_{A}= x_{A}) &=& \mathbb{P}(\phi_k(x_ {\pa_{G_{\mathcal{C}}}(k)} , \epsilon_k) \in K | X_{A}= x_{A}) \\
&=& \mathbb{P}(\phi_k(x_ {\pa_{G_{\mathcal{C}}}(k)} , \epsilon_k) \in K),
\end{eqnarray*}
where  the last equality follows from the fact the graph is ancestral, and thus $\epsilon_k$ is independent of $X_A$.

Suppose that $i\in \an_{G_{\mathcal{C}}}(k)$. Similarly, upon standard intervention, on $i$, since all the parents are given, we still have
$$\mathbb{P}(\tilde{X}_k \in K| \tilde{X}_{A}= x_{A})=\mathbb{P}(\phi_k(x_ {\pa_{G_{\mathcal{C}}}(k)} , \epsilon_k) \in K).$$

If $i \not\in \an_{G_{\mathcal{C}}}(k)$, then  upon standard intervention on $i$, we have that $(\epsilon_k,\tilde{X}_{\an_{G_{\mathcal{C}}}(k)},\tilde{X}_i )$ are independent, so that
\begin{equation*}
\mathbb{P}(\tilde{X}_k \in K| \tilde{X}_{A}=x_{A}, \tilde{X}_i=x_i)  = \mathbb{P}(\phi_k(x_ {\pa_{G_{\mathcal{C}}}(k)} , \epsilon_k) \in K).
\qedhere
\end{equation*}
\end{proof}

\begin{lemma}
\label{scm-fact-biv}
Let $\mathcal{C}$ be an SCM associated to a directed ancestral graph $G_{\mathcal{C}}$.  Consider  distinct $i,j,k \in V$ such that
$j\notin\pa_{G_{\mathcal{C}}}(k)$ and $k\notin\pa_{G_{\mathcal{C}}}(j)$,
and  $A \subseteq V$ be such that $\pa(j,k) \subseteq A \subseteq \an_{G_{\mathcal{C}}}(j,k)$.
 There exist regular conditional probabilities such that their marginals on $\mathcal{X}_{j,k}$ satisfy

$$P^{j,k}( J, K \cd x_{A}) = \mathbb{P}\big(\phi_j(x_ {\pa_{G_{\mathcal{C}}}(j)} , \epsilon_j) \in J), \phi_k(x_ {\pa_{G_{\mathcal{C}}}(k)} , \epsilon_k) \in K\big),$$

\begin{equation*}
P^{j,k}_{\doo(i)}(J,K \cd x_{A}) = P^{j,k}(J,K \cd x_A) \text{ \ for all  $i  \in \an_{G_{\mathcal{C}}}(j,k)$,    $x_A \in \mathcal{X}_A$},
\end{equation*}
and
\begin{equation*}
P^k_{\doo(i)}( J,K \cd x_{A}, x_i) = P^k(J,K  \cd x_A) \text{ \ for all  $i \not \in  \an_{G_{\mathcal{C}}}(j,k)$, $x_i \in \mathcal{X}_i$,    $x_A \in \mathcal{X}_A$},
\end{equation*}

for all $J \subseteq \mathcal{X}_j$ and $K \subseteq \mathcal{X}_k$  measurable.
\end{lemma}
\begin{proof}
We have
\begin{eqnarray*}
\mathbb{P}(X_j \in J, X_k \in K \cd X_{A}= x_{A}) &=& \mathbb{P}(\phi_j(x_ {\pa_{G_{\mathcal{C}}}(j)} , \epsilon_j) \in J, \phi_j(x_ {\pa_{G_{\mathcal{C}}}(k)} , \epsilon_k) \in K \cd X_{A}= x_{A}) \\
&=& \mathbb{P}\big(\phi_j(x_ {\pa_{G_{\mathcal{C}}}(j)} , \epsilon_j) \in J), \phi_k(x_ {\pa_{G_{\mathcal{C}}}(k)} , \epsilon_k) \in K\big),
\end{eqnarray*}
since by definition $A$ does not contain $j$ nor $k$, and thus $X_A$ and $( \epsilon_j, \epsilon_k)$ are independent.

Suppose that $i\in  \an_{G_{\mathcal{C}}}(j,k)$. Similarly, upon standard intervention, on $i$, since, we assumed
$j\notin\pa_{G_{\mathcal{C}}}(k)$ and $k\notin\pa_{G_{\mathcal{C}}}(j)$,
all the parents are given, \
we still have
$$\mathbb{P}(\tilde{X}_j \in J, \tilde{X}_k \in K \cd \tilde{X}_{A}= x_{A})=\mathbb{P}(\phi_k(x_ {\pa_{G_{\mathcal{C}}}(k)} , \epsilon_k) \in K, \phi_j(x_ {\pa_{G_{\mathcal{C}}}(j)} , \epsilon_j) \in J).$$

If $i \not\in \an_{G_{\mathcal{C}}}(j,k)$, then  upon standard intervention, on $i$, we have  $( (\epsilon_j, \epsilon_k), \tilde{X}_{\an_{G_{\mathcal{C}}}(j,k)},\tilde{X}_i )$ are independent,  so that
\begin{equation*}
\mathbb{P}(\tilde{X}_j \in J, \tilde{X}_k \in K \cd \tilde{X}_{A}=x_{A}, \tilde{X}_i=x_i)  = \mathbb{P}(\phi_k(x_ {\pa_{G_{\mathcal{C}}}(j)} , \epsilon_j) \in J, \phi_k(x_ {\pa_{G_{\mathcal{C}}}(k)} , \epsilon_k) \in K)
\qedhere
\end{equation*}
\end{proof}

\begin{theorem}\label{coro:scmrvariate}
Let $\mathcal{C}$ be an SCM associated to a directed ancestral graph $G_{\mathcal{C}}$. Assume that $\{\mathcal{P}_{\doo}(\mathcal{C}), P_{\mathcal{C}}\}$ is compatible, and $\mathcal{P}_{\doo}(\mathcal{C})$  satisfies the edge-cause condition w.r.t.\ $G_{\mathcal{C}}$. Then  $\mathcal{P}_{\doo}(\mathcal{C})$ is a bivariate-quantifiable interventional family.
\end{theorem}
\begin{proof}
Using Proposition \ref{prop:causeanscm}, we have  that Lemmas \ref{scm-fact} and \ref{scm-fact-biv} imply Axioms \ref{eq:ax1} and \ref{eq:ax2n}.
\end{proof}

\subsection{Characterizing the edge-cause condition}
\label{section-jlmr}
Consider an SCM, where $i$ is a parent of $j$.
Previously \cite[Proposition 35]{sad22}, we characterized completely which  functions $\phi_j$ would result in
the
independence $i \ci j \cd \an(i,j)$ and thus the converse pairwise Markov property.    Using similar ideas, we give a  similar characterization of when the edge-cause holds using the following two remarks.

\begin{proposition}
	\label{straight}
	 Consider an SCM, with $i \in \pa(j)$.
Write  $X_j= \phi_j(X_{\pa(j) \setminus \{i\} }, X_i, \epsilon_j)$.
        For almost all $x_i \in \mathcal{X}_i$, let $D_{x_i}$ be a random variable with law $Q_{x_i}(\cdot) = \mathbb{P}(X_{\pa(j) \setminus \{i\} } \in \cdot \cd X_i = x_i)$ that is independent of $\epsilon_j$.  Then $X_i$ is independent of $X_j$ if and only if   for almost all  $x_i \in \mathcal{X}_i$, we have that
	 $X_j$ has the same law as $\phi_j(D_{x_i}, x_i, \epsilon_j).$
	\end{proposition}
\begin{proof}
	Let $\mu_{i}$ denote the law of $X_i$ and let $\E_{\epsilon_j}$ denote the law of $\epsilon_j$. For measurable $F_i \subset \mathcal{X}_i$ and  $F_j \subset \mathcal{X}_j$, we have
	\begin{eqnarray}
		\mathbb{P}(X_i &\in& F_i, X_j \in F_j) = \nonumber \\
  &&\int\int\int \mathbf{1}[ \phi_j(x_{\pa(j) \setminus \{i\}} , x_i, e_j) \in F_j)] \mathbf{1}[x_i \in F_i]dE_{\epsilon_j}(e_j) dQ_{x_i}(x_{\pa(j) \setminus \{i\}} ) d\mu_i(x_i)   \nonumber\\
		\label{back-sub}
		&=& \int \mathbf{1}[x_i \in F_i] \mathbb{P}\big(  \phi_j(D_{x_i}, x_i, \epsilon_j) \in F_j\big) d\mu_i(x_i).
		\end{eqnarray}
	Thus if $X_j$ and $X_i$ are independent,  then
	$$\mathbb{P}(X_j \in F_j) = \int \frac{\mathbf{1}[x_i \in F_i]}{\mathbb{P}(X_i \in F_i)} \mathbb{P}\big(  \phi_j(D_{x_i}, x_i, \epsilon_j) \in F_j\big) d\mu_i(x_i),$$
	for all measurable $F_i \subset \mathcal{X}_i$ and all  $F_j \subset \mathcal{X}_j$ with nonzero probability; since this holds for \emph{all} measurable subsets $F_i$ with nonzero probability,  independence implies that, for almost every $x_i \in F_i$, we have
	\begin{equation}
				\label{dist-for-sub}
	\mathbb{P}(X_j \in F_j) = \mathbb{P}\big(  \phi_j(D_{x_i}, x_i, \epsilon_j) \in F_j\big)
	\end{equation}
	and they have same law, as desired.

	 On the other hand, substituting \eqref{dist-for-sub} in \eqref{back-sub}, we obtain the independence of $X_i$ and $X_j$.
	\end{proof}	
	
	\begin{remark}[Characterizing the edge-cause condition in $P_{\doo(i)}$: cause]
	\label{straightt}
	Notice that Proposition \ref{straight} remains valid for the associated SCM that is defined by intervention on $i$.  In particular, with the intervention $X_i = \tilde{X}_i$, we replace $Q_{x_i}$ by
	$$\tilde{Q}_{x_i}(\cdot) = \mathbb{P}(\tilde{X}_{\pa(j) \setminus \{i\} } \in \cdot \cd \tilde{X}_i = x_i).$$
	We remark that it is possible on intervention, $\tilde{X}_{\pa(j) \setminus \{i\}}$ is independent $\tilde{X}_i$.  Thus $\tilde{Q}_{x_i}$ and ${Q}_{x_i}$ could behave quite differently.  \erk
	\end{remark}

\begin{proposition}
	\label{pairwise-distinct}
	  Consider an  SCM, with $i \in \pa(j)$.
   Write  $X_j= \phi_j(X_{\pa(j) \setminus \{i\} }, X_i, \epsilon_j)$.
   %
   Suppose that $ S \supseteq  \pa(j) \setminus \ns{i}$.    Then   $i\ci j\cd S$ if and only if for almost all $x_{s} \in \mathcal{X}_{S}$ the conditional law of $X_j$ given $X_{S} = x_{S}$ is the same as the law of $\phi_j(x_{\pa(j) \setminus \{i\}}, x_i, \epsilon_j)$,
	for almost all $x_i \in \mathcal{X}_i$ in the support of the conditional distribution of $X_i$ given $x_{S}$.
	\end{proposition}

\begin{proof}
The proof is similar to that of Proposition \ref{straight} and is a routine variation of the  \cite[Proposition 35]{sad22}.
\end{proof}

\begin{remark}[Characterizing the edge-cause condition in $P_{\doo(i)}$: direct cause]
As in Remark \ref{straightt}, we apply Proposition \ref{pairwise-distinct} to the  associated SCM that is defined by intervention on $i$; we set $S =  \icause_i(j)\setminus\{i\}$. \erk
\end{remark}
%
%

\section{Identifying cases that need extra or multiple concurrent interventions}\label{sec:intmult}
Our theory is based on the interventional family $\mathcal{P}_{\doo}=\{P_{\doo(i)}\}_{i\in V}$, which only allows single interventions.

We note again that, for ancestral causal graphs, under certain conditions, direct cause can be simply defined using $i\notci_{P_{\doo(i)}} k\cd \cause{k}\setminus\{i\}$; see Section \ref{sec:anc}.  The following example shows that for the case of SCMs for non-ancestral graphs, this definition might misidentify some direct causes.
\begin{example}
    In Example \ref{ex:dcauseit} and the graph of Figure \ref{fig:dcause}, we observed that an iterative procedure is needed to obtain $i\notci_{P_{\doo(i)}} k\cd \icause_i(k)\setminus\{i\}\}$, which does not coincide with $i\notci_{P_{\doo(i)}} k\cd \cause{k}\setminus\{i\}$ in this case.

    In an SCM associated to the below graph with standard intervention, we see that even $i\notci_{P_{\doo(i)}} k\cd \icause_i(k)\setminus\{i\}$ misidentifies $i$ as the direct cause of $k$ since $i$ has no parents in the graph. We observe that, in this case the independence $i\ci_{P_{\doo(j)}} k\cd \cause{k}\setminus\{i\}$ holds.
    \begin{figure}[htb]
\centering
\begin{tikzpicture}[node distance = 8mm and 8mm, minimum width = 6mm]
    \begin{scope}
      \tikzstyle{every node} = [shape = circle,
      font = \scriptsize,
      minimum height = 6mm,
      inner sep = 2pt,
      draw = black,
      fill = white,
      anchor = center],
      text centered]
      \node(j) at (0,0) {$j$};
      \node(h) [right = of j] {$h$};
      \node(k) [right = of h] {$k$};
      \node(i) [left = of l] {$i$};
    \end{scope}
		
    \begin{scope}[<->, > = latex']
    \draw (j) to [bend left]  (k);
    \end{scope}
    \begin{scope}[->, > = latex']
    \draw (i) -- (j);
    \draw (j) -- (h);
    \draw (h) -- (k);
    \end{scope}

    \end{tikzpicture}
		\caption{{\footnotesize A (non-maximal) graph associated to an SCM}}\label{fig:nonmaxi}
		\end{figure}
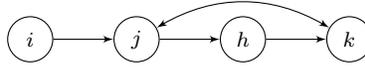
  \erk
\end{example}

Below, we classify the cases for non-ancestral graphs, where direct cause cannot be defined in the way described above.

The result below shows that all cases where $\cause{k}$ is not sufficient to define $\dcause{k}$ result in PIPs:

\begin{proposition}
For a transitive interventional family $\mathcal{P}_{\doo}$, assume $P_{\doo(i)}$ is Markovian to the $i$-intervened graph, $G_i(\mathcal{P}_{\doo})$ (e.g., by satisfying conditions of Theorem \ref{thm:markov-int}.) For non-adjacent pair of nodes $i$ and $k$, let $i\in\cause{k}\setminus\dcause{k}$. If $i\notci_{P_{\doo(i)}} k\cd \cause{k}\setminus\{i\}$, then there is a PIP between $i$ and $k$ in $G(\mathcal{P}_{\doo})$.
\end{proposition}
\begin{proof}
Since  $P_{\doo(i)}$ is Markovian to $G_i(\mathcal{P}_{\doo})$, there is a $\sigma$-connecting path between $i$ and $k$ given $\cause{k}\setminus\{i\}$;  this path is also connecting in $G(\mathcal{P}_{\doo})$. Since $\cause{k}=\an_{G(\mathcal{P}_{\doo})}(k)$, and the
global Markov property
does not imply the pairwise Markov property only for inseparable pairs, $i,k$ must be an inseparable pair. By Proposition \ref{prop:nonmaxp}, there is a PIP in $G(\mathcal{P}_{\doo})$ between $i$ and $k$.
\end{proof}

For $i\in\cause{k}$, there could be two types of  PIPs between $i$ and $k$, described in the above proposition.    Notice again that such cases can happen only for non-ancestral graphs and when the true causal graph is non-maximal.
\begin{enumerate}[(1) ]
\item
If this PIP is not a PIP in $G_i(\mathcal{P}_{\doo})$ (as in the case of Figure \ref{fig:dcause}), i.e., an inner node of the PIP is only an ancestor of $k$ via $i$, then the iterative procedure to define direct cause works as it is designed to ensure that a direct cause is not placed in the causal graph incorrectly.
\item
If the PIP is a PIP in $G_i(\mathcal{P}_{\doo})$ (as is, e.g., the case in Figure \ref{fig:nonmaxi}), then the current theory is incomplete in the sense that some direct causes may be misidentified as our theory considers $i\in\dcause{k}$, and it always places an arrow from $i$ to $k$ in $G(\mathcal{P}_{\doo})$. We have two sub-cases here:

\begin{enumerate}[(a) ]
\item If there exists maximum one such PIP  between each pair of nodes, then we propose the following adjustment to the definition of direct cause, which fixes this issue:

We define $i$ to be a direct cause of $k$ if for every $j$ (that may be $i$ but not $k$), it holds that $i\notci_{P_{\doo(j)}} k\cd \cause{k}\setminus\{i\}$.

Hence, as a procedure, one can generate the intervened causal structure $S_i(\mathcal{P}_{\doo})$, and if there is a PIP between $i$ and $k$ then performs the extra test $$i\ci_{P_{\doo(j)}} k\cd \cause{k}\setminus\{i\},$$
for a $j$ on PIP. If it holds, then the arrow from $i$ to $k$ will be removed.

We note that, for observable interventional families, and for separable pairs, by Axioms \ref{eq:ax20} and \ref{eq:ax2n},  the original definition of direct cause is equivalent to this new one. However, this definition, is an extension of the original definition for inseparable pairs.

Notice also that in such cases, similar to bows, by only knowing $\mathcal{P}_{\doo}$, it is not possible to distinguish an arrow from an arc between $i$ and $k$. We treat such cases as a direct cause from $i$ to $k$.

\item  If there are more than one such PIPs between $i$ and $k$ then, In the case of faithful SCMs, no matter which $j$ one intervenes on, $i\notdse_{\sigma} k\cd \cause{k}\setminus\{i\}$. Hence, Markov property does not imply independence of $i$ and $k$ given $\cause{k}\setminus\{i\}$.

In such cases, one should intervene concurrently on one node on each of these PIPs to determine whether $i$ and $k$ become separated given $\cause{k}\setminus\{i\}$. Hence our theory based on single interventions is not identifying such direct causes. To be more precise, we have the following remark.
\end{enumerate}
\end{enumerate}

\begin{remark}
Let $\mathcal{P}_{\doo}$ be a transitive interventional family. Assume $i\in\cause{k}$. Suppose that there are $r$ PIPs $\pi_1,\cdots,\pi_r$ between $i$ and $k$. In order to identify whether $i\in\dcause{k}$, the existence of a concurrent intervention $P_{\doo(j_1,\cdots,j_r)}$ of size $r$, is necessary, where each $j_s$, $1\leq s\leq r$, is an inner node of $\pi_s$. \erk
\end{remark}

\section{Summary and discussion}\label{sec:summary}
For a given family of distributions for a finite number of variables with the same state space (which we refer to as interventional family), we have defined \emph{the} set of causes of each variable. This is equivalent to the definition of the cause provided in the literature for SCMs. We have provided transitivity of the cause as the first axiom on the given family, as it leads to reasonable causal graphs.   We have provided weaker conditions for interventional families under which the family is transitive if the interventional distributions satisfy the property of singleton-transitivity.  Although, in general, causation does not seem to be transitive \cite{hal00}, it seems to us that examples for which causation is not transitive mostly follow the counterfactual interpretation of causation rather than the conditional independence one; and possible examples based on conditional independence interpretations do not satisfy singleton-transitivity. We refrain from philosophical discussions here, but note that in our opinion, representing causal structure using graphs implicitly implies that one is focusing on the cases where cause in indeed transitive as the directed paths in directed graphs are transitive.

We have defined the concept of direct cause using the definition of the cause and pairwise independencies given the joint causes in $P_{\doo(i)}$ with an iterative procedure in the general case. However, under the assumption that the causal graph is directed ancestral, causal graphs can be defined without an iterative procedure. The major departure here from the literature is that the direct cause of a variable is defined using single interventions and by \emph{conditioning} on other causes of the variable. As opposed to the defined causal relationships, which stand true in a larger system of variables, the direct causal relationship clearly depends on the system\textemdash it seems that one can always add a new variable to the system that breaks the direct causal relationship by sitting between the two variables as an intermediary.

Our original motivation to write this paper was to relax the common assumption that there \emph{exists} a true causal graph, and, thereby, the goal of causal inference is solely to learn or estimate this graph. We do not need any such assumption under our axiomatization as we define causal graphs using intervened graphs, which themselves are defined using the concept of the direct cause for arrows and the pairwise dependencies given the joint causes in the interventional distributions for arcs. Arcs represent latent variables, and our generated graphs also allow for causal cycles. Transitivity ensures that the causal relationships in the interventional family and the ancestral relationships in graph are interchangeable.

We believe this setting can be extended to causal graphs that unify \emph{anterial graphs} \cite{sadl16} with cyclic graphs. Such a graph represents, in addition, symmetric causal relationships implied by feedback loops; see \cite{lau02}. In order to do so, some extension of Markov properties, presented here for BDMGs, for this larger class of graphs is needed.


For the case where causal cycles exist, one advantage of the setting presented here is that it is easy to provide examples for cyclic graphs under our axiomatization; see Example \ref{all-example}.   This is in contrast with the case of SCMs with cycles, where, for this purpose, strong solvability assumptions must be satisfied \cite{bon20}.

We notice again that there is no need for $P$ to define the causal graph. We provide a  minimalist and a maximalist approach to place an arc in the casual graph based on whether the arc exists in intervened graphs. This will only lead to different causal graphs where arcs between inseparable pairs of nodes are or are not present.

We have also axiomatized interventional families with an observed distribution relating certain conditional independence properties in $P$ and $P_{\doo(i)}$, and called a family that satisfies these axioms (strongly) observable interventional families. Although we have defined the causal graph only using $P_{\doo(i)}$, we show that under these axioms, the arcs in causal graph can be directly defined using the observed distribution $P$.


There are two main results in this paper: firstly, under compositional graphoids with transitive  families, any  interventional distribution  $P_{\doo(i)}$ is Markovian to the intervened causal graph over node $i$; secondly, under observable interventional families, $P$ is Markovian to the generated causal graph. We note that our definition of the arc, partly ensures automatically that the Markov properties are satisfied. However, under an alternative definition of the arc we provided\textemdash which places an arc when the endpoint variables are always dependent regardless of what we condition on\textemdash no assumptions related to Markov property is in place\textemdash in this case, we need the extra assumption of ordered upward- and downward-stability to obtain the Markov property. 
 These Markov-property results are analogous to the case of SCMs, and, consequently, this allows the developed theory for causality in the SCMs to be embedded in the general setting of this paper.

We mostly work on definitions and axioms for interventional families that are only related to conditional independence structure of interventional and observed distributions.  We note that although they are sufficient for generating and making sense of causal graphs, for ``measuring'' causal effects (which we do not discuss) they are not sufficient. For that purpose, for the case of directed ancestral causal graphs, we provide the axiom of (bivariate) quantifiable interventional families, which relates the univariate (and bivariate) conditional-marginals of the distributions in the family to those of an underlying distribution $P$. In principle, $P$ could be learned via observation, and in the case of DAGs, it is determined uniquely by the interventional family. The extension of this axiom to BDMGs seems quite technical, and requires further study.

Note that satisfaction of the axioms for a family of distributions does not mean that the family provides the \emph{correct} interventions\textemdash refer again to Example \ref{all-example} to observe that all three types of edges, as the causal graph for different interventional families of two variables with the same underlying distribution $P$, can occur. Finding the correct interventional families is a question for mathematical and statistical modeling.  One can think of this as being analogous to Kolmogorov probability axioms \cite{kol60}: a
measurable space
satisfying Kolmogorov axioms does not mean that it provides the correct probability for the experiment at hand. This is not the case in the SCM setting, as in the presence of densities,  interventional distributions with full support are equivalent \cite{pet17}\textemdash this is because the causal graph in this setting is assumed to exist and already set in place. Example \ref{two-triangles} shows that even the skeleton of the causal graph could change by the change of interventions with the same underlying distribution.

When we relate SCMs to the setting of this paper, we find that if the SCM satisfies some weaker version of faithfulness given by the edge-cause and converse pairwise Markov property, then, in the case where the natural graph associated with the SCM is ancestral, the causal graph, given by standard interventions, is the same as the SCM graph; if the graph is not ancestral, then if the interventions are transitive, we can recover this result for maximal BDMGs.  These results  demonstrate that our  theory is compatible with the standard theory, for a large class of SCM.  

We have not provided conditions on SCMs under which the cause is transitive, although it is not used for the main results related to SCMs being embedded in the setting of this paper. Our initial investigation revealed that this is quite a technical problem. This is nevertheless beyond the scope of this paper, and is a subject of future work.

Finally, although an advantage of our theory is that it only relies on single interventions, our theory might misidentify  direct causes between primitive inducing paths in the intervened graphs for non-maximal non-ancestral causal graphs. We provide an adjustment to deal with this when there is only one PIP exists between a pair. If there are more PIPs between a pair, we showed that we need multiple concurrent intervention of the size of the number of PIPs between the pair.

Similarly, our theory does not include some cases where multiple concurrent interventions could act as the cause of a random variable whereas none of them individually act as the cause; for example; see Example \ref{ex:doubledo}. Understanding these cases, and developing a similar theory for such cases is a subject of future work.

\section{Appendix}
\label{appn}
\section*{Pairwise Markov properties for bowless directed mixed graphs (BDMGs)} 
Here, we prove the equivalence of the pairwise Markov property \eqref{P} and the global Markov property by defining an intermediary pairwise Markov property for
acyclic directed mixed graphs
(ADMGs).

First, we need to define the concept of acyclification from \cite{bon20,forre2017markov}, which generates a Markov-equivalent acyclic graph from a graph that contains directed cycles by adding and replacing some edges in the original graph:  For a directed mixed graph $G$ (which contains directed cycles), the \emph{acyclification} of $G$  is the acyclic graph $G^{\mathrm{acy}}$ with the same node set, and the following edge set:
\begin{enumerate}
\item[$\bullet$  \ ]
There is an arrow $j\fra i$ in $G^{\mathrm{acy}}$ if $j\in\pa(\sco(i))\setminus\sco(i)$ in $G$;
\item[$\bullet$  \ ]
and there is an arc $i\arc j$ in $G^{\mathrm{acy}}$ if and only if there exist $i'\in\sco(i)$ and $j'\in\sco(j)$ with $i'=j'$ or $i'\arc j'$ in $G$.
\end{enumerate}
\begin{lemma}
\label{lem:acy}
Let $G$ be a BDMG.   For distinct nodes $i$ and $j$, we have  $i\in\an_G(j)$ and $i\notin\sco_G(j)$  if and only if  $i\in\an_{G^{\mathrm{acy}}}(j)$.
\end{lemma}
\begin{proof}
($\Rightarrow$) Suppose that there is a directed path $\pi=\langle i=i_0,i_1\cdots,i_n=j\rangle$ from $i$ to $j$ in $G$. This path is a directed path in $G^{\mathrm{acy}}$ unless an arrow $i_ri_{r+1}$ turns into an arc. However, in this case, $i_{r+1}\in\sco(i_r)$ and, hence,  $i_{r-1}\in\pa(i_{r+1})$. An inductive argument implies the result since $i\notin\sco(j)$.

($\Leftarrow$) Conversely, if $\pi=\langle i=i_0,i_1\cdots,i_n=j\rangle$ is a directed path in $G^{\mathrm{acy}}$ from $i$ to $j$ then an arrow $i_ri_{r+1}$ on $\pi$ exists in $G$ unless $i_r$ is a parent of another node $k$ such that $k\in\sco(i_{r+1})$. This implies that $i_r\in\an(i_{r+1})$ in $G$, and, again, an inductive argument implies that $i\in\an(j)\setminus\pa(j)$.

In addition, $i\notin\sco(j)$ in $G$ since, if, for contradiction, this is not the case then all the arrows in the strongly connected component containing $i$ and $j$ turn into arcs. If there is an arrow $k\ell$ generated in $G^{\mathrm{acy}}$ that makes the directed path from $i$ to $j$,  then by the construction of acyclification, it can be seen that $k$ is an ancestor of $\ell$, and hence they are still a part of the same strongly connected component, and they must turn into arcs, which is a contradiction.
\end{proof}
\begin{proposition}
\label{prop:anc}
If $G$ is a BDMG, then $G^{\mathrm{acy}}$ is an ADMG.
\end{proposition}
\begin{proof}
It is easy to see from the definition of acyclification that  $G^{\mathrm{acy}}$ is acyclic; see also \cite{bon20}.
\end{proof}

It was shown in \cite[Proposition A.19]{bon20} that the global Markov property could be read off equivalently from the acyclification of a directed mixed graph; see also \cite{forre2017markov}.
\begin{proposition}[Equivalence of the separation criteria \cite{bon20}]
\label{prop:sepacy}
Let $G$ be a BDMG. For disjoint  subsets of nodes $A,B,C$, we have
$$A\dse_{\sigma}B\cd C \text{ in } G \iff A\dse_m B\cd C \text{ in } G^{\mathrm{acy}}.$$
\end{proposition}

Before proceeding to provide required results for proving Theorem \ref{prop:pairMark}, we prove Proposition \ref{prop:nonmaxp}. It was proven in \cite{sadl16} that every inseparable pair is connected by a PIP for a generalization of ADMGs; see Section 4 of this paper. The definition of a PIP for ADMGs is as follows: a path is a PIP if every inner node is a collider node and an ancestor of one of the endpoints.

\begin{proof}[Proof of Proposition \ref{prop:nonmaxp}]
We suppose that $(i,k)$ is an inseparable pair in $G$. In $G^{\mathrm{acy}}$, by Propositions \ref{prop:anc} and \ref{prop:sepacy}, this is an inseparable pair or there is an edge between $i$ and $k$.

First we show that if there is an edge between $h$ and $\ell$ in $G^{\mathrm{acy}}$ then this is an edge in $G$ or there is a PIP in $G$. Assume it is not an edge in $G$. There are two cases: If it is an arrow, then $h$ and $\ell$ are part of the same strongly connected component and a path on this is a PIP in $G$ by definition. If $h\ell$ is an arc, then $h$ and $\ell$ are part of the same strongly connected component again or $h,h'$ and $\ell,\ell'$  are parts of the same strongly connected component, respectively, where there is an arc between $h'$ and $\ell'$. In this case, a path consisting of a directed path on the component between $h$ and $h'$, the $h'\ell'$-arc, and a directed path on the component between $\ell$ and $\ell'$ is a PIP in $G$.

The above will resolve the case where there is an edge between $i$ and $k$. Now consider the case where  $i$ and $k$ are inseparable in $G^{\mathrm{acy}}$. It was proven in \cite{sadl16} that every inseparable pair is connected by a PIP for a generalization of ADMGs; see Section 4 of this paper.

Hence, we need to show that if there is a primitive inducing path $\pi$ in $G^{\mathrm{acy}}$ between $i$ and $k$ then there is a PIP in $G$ between $i$ and $k$.
From the previous argument about edges in $G^{\mathrm{acy}}$,  there are edges or PIPs corresponding to every edge of $\pi$. We put all these  together to form a walk between $i$ and $k$. We then look at a subpath of this walk between $i$ and $k$, and call it $\rho$.

Every edge on $\rho$ is either an arrow, where endpoints are in the same strongly connected component or an arc as required.

Hence, we only need to show that every inner node is in $\an_G(\{i,k\})$. We know that they are all in $\an_{G^{\mathrm{acy}}}(\{i,k\})$;
thus, the result follows from
Lemma \ref{lem:acy}.
%
%
\end{proof}
We now recall that a pairwise Markov property for ADMGs is the same as in the BDMGs as defined by  \eqref{P}.
%
%
%
The proof of the equivalence of the pairwise Markov property and the global Markov property for ADMGs is a trivial special case of \cite[Corollary 5]{sadl16}.
\begin{proposition}[Equivalence of pairwise and global Markov properties for ancestral graphs \cite{sadl16}]
\label{prop:pairsMarkeq}
Let $G$ be an ADMG, and assume that the distribution $P$ satisfies the intersection property and the composition property.  Then $P$ satisfies the pairwise Markov property  w.r.t.\ $G$  if and only if it is Markovian to $G$.
\end{proposition}
\begin{proposition}
\label{prop:pairgloMarkeq}
Let $G$ be a BDMG.  
If $P$ satisfies the pairwise Markov property w.r.t.\ $G$,  then it satisfies the pairwise Markov property w.r.t.\ its acyclification,  $G^{\mathrm{acy}}$.
\end{proposition}
\begin{proof}
Consider two arbitrary non-adjacent nodes $i$ and $j$ in $G^{\mathrm{acy}}$; they are not adjacent in $G$, and also $i\notin\sco_G(j)$ since otherwise they would have made adjacent after acyclification. Thus,  by the pairwise Markov property for $G$, we have
$$i\ci_P j\cd [\an_G(i)\cup\an_G(j)]\setminus\{i,j\}.$$
Moreover, by Lemma \ref{lem:acy}, $\an_G(i)=\an_{G^{\mathrm{acy}}}(i)$, and same equality holds for the node $j$.  Thus, we obtain the pairwise Markov property for $G^{\mathrm{acy}}$.
\end{proof}
\begin{proof}[Proof of Theorem \ref{prop:pairMark}]
 By Proposition \ref{prop:pairgloMarkeq}, the pairwise Markov property w.r.t.\ $G$ carries over to its acyclification $G^{\mathrm{acy}}$, which by Proposition \ref{prop:pairsMarkeq} implies the global Markov property  w.r.t.\ $G^{\mathrm{acy}}$, since by Proposition \ref{prop:anc}, $G^{\mathrm{acy}}$ is an ADMG.  Finally, by  Proposition \ref{prop:sepacy}, we recover the global Markov property w.r.t.\ the original graph $G$.
\end{proof}
\begin{proof}[Proof of Proposition \ref{prop:Markpair}]
It is enough to prove, for two arbitrary nodes $i$ and $j$ in a bowless directed mixed graph $G$, that $i\dse_{\sigma}j\cd \an(\{i,j\})$. By Proposition \ref{prop:sepacy}, this is equivalent to $i\dse_m j\cd \an(\{i,j\})$ in $G^{\mathrm{acy}}$. Notice that, by Proposition \ref{prop:anc}, $G^{\mathrm{acy}}$ is an ADMG. In addition, by Proposition \ref{prop:sepacy}, for every separation statement between two nodes $k$ and $l$ in $G$, there is a separation statement between $k$ and $l$ in $G^{\mathrm{acy}}$; hence, maximality of $G$ implies that $G^{\mathrm{acy}}$ is maximal. The result now follows from the fact that, for maximal ADMGs, the separation $i\dse_m j\cd \an(\{i,j\})$ always holds; see \cite{sadl16}.  
\end{proof}

\section{Acknowledgments}
The authors are grateful to Patrick Forr\'{e} for a helpful discussion on pairwise Markov properties for graphs with directed cycles, and to Philip Dawid, Thomas Richardson, and Jiji Zhang for helpful general discussions related to this manuscript. Work of the first author is supported in part by the EPSRC grant EP/W015684/1.

\bibliographystyle{imsart-number} 
\bibliography{bib}       

\end{document}